\documentclass[12pt]{article}
\usepackage[margin=1in, a4paper]{geometry}

\usepackage{amsmath,enumerate,amsfonts,color,amssymb,amsthm,txfonts,cancel}

\usepackage[normalem]{ulem}

\usepackage{hyperref}

\usepackage{tikz}

\usepackage[framemethod=tikz]{mdframed}

\def\RR{{\mathbb R}}

\def\mcF{{\mycal F}}

\def\ds{\displaystyle}

\def\eps{{\varepsilon}}
\def\vareps{{\varepsilon}}
\def\loc{{\rm loc}}
\def\mcF{{\mycal F}}
\newtheorem{theorem} {\sc  Theorem\rm} [section]

\newtheorem{corollary} [theorem] {\sc  Corollary\rm}
\newtheorem{lemma} [theorem] {\sc  Lemma\rm}
\newtheorem{proposition} [theorem] {\sc  Proposition\rm}
\newtheorem{prop} [theorem] {\sc  Proposition\rm}

\newtheorem{remark}[theorem]{\sc  Remark\rm}

\def\bproof{\noindent{\bf Proof.\;}}
\def\eproof{\hfill$\square$\medskip}

\def\nd{\noindent}

\newcounter{marnote}

\newcommand{\norm}[1]{\left\Vert#1\right\Vert}

\DeclareFontFamily{OT1}{rsfs}{}
\DeclareFontShape{OT1}{rsfs}{m}{n}{ <-7> rsfs5 <7-10> rsfs7 <10-> rsfs10}{}
\DeclareMathAlphabet{\mycal}{OT1}{rsfs}{m}{n}

\def\dist{{\rm dist}}

\def\supp{{\rm Supp}\,}

\def\be{\begin{equation}}
\def\ee{\end{equation}}

\def\dist{{\rm dist}}

\def\supp{{\rm Supp}\,}

\newcommand{\R}{\mathbb{R}}
\newcommand{\Ss}{\mathbb{S}}

\newcommand{\e}{\varepsilon}

\newcommand{\abs}[1]{\left\vert{#1}\right\vert}
\def\be{\begin{equation}}
\def\ee{\end{equation}}
\def\bea#1\eea{\begin{align}#1\end{align}}

\numberwithin{equation}{section}

\mdfdefinestyle{MyFrame}{%
	linecolor=black,
	outerlinewidth=1pt,
	roundcorner=5pt,
	innertopmargin=\baselineskip,
	innerbottommargin=\baselineskip,
	innerrightmargin=10pt,
	innerleftmargin=10pt,
	backgroundcolor=white!20!white}

\begin{document}
\title{The Ginzburg-Landau system with general potential: maximum principle and gradient estimates}

\author{Radu Ignat\thanks{Institut de Math\'ematiques de Toulouse, UMR 5219, Universit\'e de Toulouse, CNRS, UPS IMT, F-31062 Toulouse Cedex 9, France.
Email: Radu.Ignat@math.univ-toulouse.fr}~, Luc Nguyen\thanks{Mathematical Institute and St Edmund Hall, University of Oxford, Andrew Wiles Building, Radcliffe Observatory Quarter, Woodstock Road, Oxford OX2 6GG, United Kingdom. Email: luc.nguyen@maths.ox.ac.uk}~, Valeriy Slastikov\thanks{School of Mathematics, University of Bristol, University Walk, Bristol, BS8 1TW, United Kingdom. Email: Valeriy.Slastikov@bristol.ac.uk}~ and Arghir Zarnescu\thanks{IKERBASQUE, Basque Foundation for Science, Plaza Euskadi 5, 48009 Bilbao, Bizkaia, Spain.}\,  \thanks{BCAM,  Basque  Center  for  Applied  Mathematics,  Mazarredo  14,  E48009  Bilbao,  Bizkaia,  Spain.
(azarnescu@bcamath.org)}\, \thanks{``Simion Stoilow" Institute of the Romanian Academy, 21 Calea Grivi\c{t}ei, 010702 Bucharest, Romania.}}

\date{}

\maketitle

\begin{abstract}
We study critical points of the Ginzburg-Landau energy functional  $$\mcF_\vareps[u] = \int_\Omega \Big[ \frac{1}{2}|\nabla u|^2 + \frac{1}{2\vareps^2} W(1 - |u|^2)\Big]\,dx, \quad u \in H^1(\Omega, \R^N),$$ 
with $\Omega \subset \R^M$, $\vareps>0$, $M,N \geq 2$ and general conditions on the non-negative potential $W$ allowing for super-quadratic behaviour near its zero set. Under a Dirichlet boundary data of unit-length on $\partial \Omega$, we prove the following maximum principle: every critical point $u_\vareps$ satisfies the global uniform bound $|u_\eps| \leq 1$ in $\Omega$. Furthermore, if a family of critical points $(u_\eps)$ converges (in energy) to a smooth $\Ss^{N-1}$-valued harmonic map in the limit $\eps \to 0$, then we prove global uniform bounds for $(\Delta u_\eps)_{\eps>0}$ in $\Omega$ and, in particular, global H\"older convergence of the gradients $(\nabla u_\eps)$ in $\Omega$ as $\eps \to 0$.  
\end{abstract}

\setcounter{tocdepth}{1}
\tableofcontents

\section{Introduction}

Let $M, N\geq 2$ be integers. For a $C^2$ bounded  domain $\Omega \subset \RR^M$ and for $\eps>0$, we consider the energy functional:
\be\label{def:mcFeps}
\mcF_\vareps[u] = \int_\Omega \Big[ \frac{1}{2}|\nabla u|^2 + \frac{1}{2\vareps^2} W(1 - |u|^2)\Big]\,dx, \quad
u \in H^1(\Omega,\RR^N),
\ee
where $W \in C^2((-\infty, 1])$ is a non-negative potential function (not necessarily convex) and the modulus $|u|$ is the Euclidean norm of $u$. We consider a family of critical points $(u_\eps)_{\eps\downarrow 0}$ of $\mcF_\eps$ in $H^1\cap L^\infty(\Omega,\RR^N)$ such that, for some limit $u_0 \in H^1(\Omega, \RR^N)$,
\be
\label{ass:convH1}
u_{\eps}\rightharpoonup u_0\textrm{ weakly in }H^1(\Omega,\RR^N) \quad \textrm{and} \quad \mcF_\eps[u_\eps]\to\frac 12\int_\Omega |\nabla u_0(x)|^2\,dx \quad \textrm{ as } \eps\to 0.
\ee 
This setup is natural in the theory of $\Gamma$-convergence. Typically as $\eps \to 0$ one derives a limiting problem for $u_0$ and then obtains insights about behavior of local minimizers (or,  more generally, of  critical points) of the $\eps$-dependent functional based on the corresponding information on local minimizers (respectively, on critical points) $u_0$ of the limiting problem,  e.g. see \cite{AlbertiBaldoOrlandi2005, JerrardSternberg2009, KohnSternberg1989,  Modica1987, SandierSerfaty2004}. Conversely, this method is also fruitful in  obtaining information on the limiting problem ($\eps=0$) using the $\eps$-dependent problem. Indeed, such approach has appeared in the construction of  minimal surfaces of co-dimension one using the (scalar) Allen-Cahn model \cite{guaraco2018min,HutToneg2000,PaRi03,to2005stable,to2012stable,Wic14} or higher co-dimension minimal surfaces using Ginzburg-Landau models \cite{BaPicod2, PhPg24}.   We emphasize that in the current paper we always consider general critical points $u_\vareps$ of $\mcF_\vareps$ which may not be local minimizers.

In the case when \eqref{ass:convH1} holds, it  is expected that under suitable conditions on the potential $W$ (see assumptions \eqref{ass:Wstructure}, \eqref{ass:Wnondeg} and \eqref{ass:convex} below) the limit  $u_0$ is a critical point of the limiting problem with $\eps=0$, that is a stationary $\Ss^{N-1}$-valued harmonic map, see Remark~\ref{Rem:u0SHM} and Appendix \ref{App:HM}. Provided that such a critical point $u_0$ is ``non-degenerate", an interesting question concerns existence and uniqueness of critical points $u_\eps$ of the $\eps$-dependent Ginzburg-Landau energy in a suitable neighborhood of $u_0$, which implicitly raises the question of the suitable topology to work on. The topological aspect of this question is precisely the focus of the current paper under the additional assumption of a smooth limit $u_0$, while the issue of existence and uniqueness is addressed in our next paper \cite{INSZ26-2nd}. 
 
In the following, we assume that $u_\eps$ satisfies the Dirichlet boundary condition\footnote{In fact, our results hold true under the boundary assumption that $(u_\eps)_{\eps>0}$ is uniformly bounded in $C^{1, \beta}(\partial \Omega, \Ss^{N-1})$.}:
$$
\quad u_\eps\big|_{\partial\Omega} =u_{bd} \in H^{1/2}(\partial \Omega, \Ss^{N-1}).
$$
A critical point $u_\eps$ of $\mcF_\eps$ with this boundary condition solves
\be
\label{eq:eps}
-\Delta u_\eps=\frac1{\eps^2}W'(1-|u_\eps|^2)u_\eps \textrm{ in } \Omega, \quad u_\eps= u_{bd} \textrm{ on } \partial\Omega.
\ee
We study uniform gradient and Laplacian estimates for $u_\eps$ as $\eps \rightarrow 0$. 
The main assumptions on the $C^2$ potential function $W$ are the following:
	 \be
\label{ass:Wstructure} 
W(0) = 0, \quad W(t) \geq 0 \text{ for } t \in (-\infty,0), \quad  W(t) > 0\textrm{ for every }t\in (0, 1], 
\ee 
and there exist $\alpha \geq 1, c_0 > 0, \kappa > 0$ so that
\begin{align}
\label{ass:Wnondeg}
&\lim_{t \rightarrow 0} \frac{W'(t) - c_0 (\alpha+1) t|t|^{\alpha-1}}{|t|^{\alpha}}   = 0
\\
\label{ass:convex}
&W''(t)\ge 0, \quad \forall |t|\leq \kappa.
\end{align}

Roughly speaking \eqref{ass:Wstructure} says that $t = 0$ is a global minimum of $W$, and \eqref{ass:Wnondeg} and \eqref{ass:convex} impose a local non-degenerate property of this minimum. 
Note that \eqref{ass:Wstructure} and \eqref{ass:convex} yield
\be\label{est:Wincreas}
W'(0)=0, W'(t)\ge 0 \text{ and }  W'(-t)\leq 0, \quad \forall t\in [0,\kappa],
\ee
and $W$ is convex near $t=0$; however, $W$ is not necessarily convex away from the origin.
An example for such potential is $W(t)=c_0|t|^{\alpha+1}$ for $t\leq 1$; the particular case $W(t) = \frac{1}{2}t^2$ for $\alpha = 1$ is of special physical importance and has been considered widely in the literature. However, this potential $W(t)=c_0|t|^{\alpha+1}$ for $\alpha \geq 1$  can only be used to model second-order phase transition. The general assumptions \eqref{ass:Wstructure}, \eqref{ass:Wnondeg}, \eqref{ass:convex} considered in the present paper are more flexible and therefore can accommodate first-order phase transitions arising in a number of physical contexts e.g. in Landau-de Gennes' theory of liquid crystals \cite[Chapter 10]{smallelastic}, or Blume-Capel's model in ferromagnetics \cite{Blume66, Capel66}. Moreover, it raises interesting mathematical questions that have been studied in related papers, for sub-quadratic potentials in the scalar case in \cite{caffarelli1995uniform} and in the vectorial case in \cite{alikakos2024, alikakos2022asymptotic} while the 
 scalar setting for super-quadratic potentials has been studied  in  \cite{dipierro2025density, dipierro2018density, savin2025}.

Following the seminal works of B\'ethuel, Brezis and H\'elein \cite{BBH93, vortices}, there have been several analytical studies on the behaviour of critical points of $\mcF_\eps$ and related energy functionals in the limit $\eps \rightarrow 0$ (see e.g. \cite{BBO01, CLR18, INSZ_CRAS18, INSZ-ENS, LinRiv2, MZ10, NZ13, Struwe94, YeZhou96} and the references therein). Although a large portion of these works remains applicable in our present context, they do not directly yield estimates which are compatible with the setting in our companion paper \cite{INSZ26-2nd} where existence and uniqueness is obtained. This is because certain aspects of these works are restricted to the case where $u_\eps$ is energy minimizing and/or the case $W$ is the standard Ginzburg-Landau potential $W(t) = \frac{1}{2}t^2$. In the generality we are considering, a couple of new issues arise which we will subsequently discuss.

An important estimate for \eqref{eq:eps} that has been used extensively in many contexts is the modulus estimate $|u_\eps| \leq 1$ that we refer here as the ``maximum principle" for the (vectorial) critical point $u_\eps$. This is easy to establish if one has that $W'(t) \leq 0$ for all $t \leq 0$. Since we assume only a local behaviour of $W$ near $t =0$ and the non-negativity of $W$ away from $t = 0$, one may contemplate a counter-example: For certain exotic potential $W$ that vanishes at some $t < 0$, one may perform a mountain-pass construction to obtain a critical point of $\mcF_\eps$ whose modulus exceeds $1$ somewhere. The first result of our paper is a new maximum modulus principle which shows that such counter-example is simply impossible. We highlight that, for this result, it suffices to assume only
\be
\label{ass:Wsweak} W(0) = 0 \text{ and }  W(t) \geq 0\textrm{ for every }t\in (-\infty, 1],
\ee
which is weaker than \eqref{ass:Wstructure}. Note here a subtle difference between \eqref{ass:Wsweak} and \eqref{ass:Wstructure}: For example, in the degenerate case $W \equiv 0$ (which satisfies \eqref{ass:Wsweak} but not \eqref{ass:Wstructure}), any constant vector is a global minimizer of $\mcF_\eps$. On the other hand, when $W$ satisfies \eqref{ass:Wstructure}, any such constant global mininizer must have modulus at least $1$.

\begin{theorem}
\label{thm:MPI}
Let $W \in C^2((-\infty,1])$ satisfy \eqref{ass:Wsweak}, $\vareps > 0$ and $\Omega\subset \R^M$ be a Lipschitz\footnote{An exterior cone condition suffices in this theorem.} bounded domain. Then every critical point $u_\eps \in H^1\cap L^\infty(\Omega, \R^N)$ of $\mcF_\eps$ such that $u_\eps \big|_{\partial \Omega} \in C^{0}(\partial \Omega, \Ss^{N-1})$ belongs to $C^2(\Omega)\cap C^0(\bar \Omega)$ and satisfies $|u_\eps|\leq 1$ in $\Omega$.
\end{theorem}

We note that \eqref{ass:Wsweak} cannot be dropped in Theorem \ref{thm:MPI}. See Proposition \ref{Prop:WsweakCE}.

Theorem \ref{thm:MPI} has physical significance. In the energy functional \eqref{def:mcFeps}, the expression for the potential energy $W$ is usually obtained by Taylor expanding\footnote{This is also known as Landau expansion.}  the true physical potential energy around the isotropic state $u = 0$, suitably truncated at a certain order so as to capture certain desired physical or phenomenological properties. As such, its expression for large $|u|$ can deviate substantially from the true expression. A natural question arises whether this approximation of the potential energy has any consequence(s) in the study of the total energy functional. Our theorem shows that, as far as critical points are concerned, this inaccuracy in the formula for $W$ for large values of $|u|$ is irrelevant. What truly matters is the behaviour of the potential energy for $u$ varying in a range containing the isotropic state $u = 0$ and the energetically preferred states $|u| = 1$ yielded by the boundary condition. 

The proof of Theorem \ref{thm:MPI} is a maximum principle argument in which we construct a family of entire solutions of the scalar PDE:
\begin{equation}
-\Delta g = \frac{1}{\eps^2} W'(1 - g^2)g \quad \text{ in } \RR^M,
	\label{Eq:SGL}
\end{equation}
and use them as barriers to deduce bounds for $|u_\eps|$. See Section \ref{Sec:MP}.

\begin{remark}\label{Rem:u0SHM}
Under \eqref{ass:convH1} and \eqref{ass:Wsweak}, it is standard to show that $u_\eps \to u_0 \text{ in } H^1(\Omega, \RR^N)$ and $W(1 - |u_0|^2) = 0$ a.e. in $\Omega$. Moreover, under \eqref{ass:Wstructure} and the boundary assumption on $u_\eps\big|_{\partial \Omega}$ in Theorem \ref{thm:MPI}, it follows that $|u_0| = 1$ a.e. in $\Omega$. In fact, $u_0$ is a stationary\footnote{ In the particular case when $u_\vareps$ is minimizing for $\mcF_\vareps$ with respect to its boundary value, then $u_0$ is a minimizing harmonic map.} $\mathbb{S}^{N-1}$-valued harmonic map (see Lemma \ref{Lem:u0ConvSHM}).  
\end{remark}

The second result of the present paper is a uniform estimate for $\Delta u_\eps$ for every $\eps>0$, assuming that the limiting harmonic map $u_0$ is regular. In particular, this result lays the $C^{1, \beta}$ topological  framework for our existence and uniqueness result in our next paper \cite{INSZ26-2nd}.
\begin{theorem} 
\label{thm:MDI}
Let $W \in C^2((-\infty,1])$ be such that \eqref{ass:Wstructure}, \eqref{ass:Wnondeg} and \eqref{ass:convex} hold, $\Omega\subset \R^M$ be a $C^2$ bounded domain and 
$u_{bd}\in C^{1, \beta}(\partial \Omega, \Ss^{N-1})$ for some $\beta\in (0,1)$.
Let $(u_\eps)_{\eps>0}\subset H^1\cap L^\infty(\Omega, \R^N)$ be a family of solutions to \eqref{eq:eps} that satisfies \eqref{ass:convH1} for a stationary harmonic map $u_0\in H^1(\Omega, \Ss^{N-1})\cap C^0(\bar \Omega, \Ss^{N-1})$ as $\eps\to 0$. Then there exists a constant $C > 0$ depending only on $M, N, W, \Omega, u_{bd}$ and $u_0$ such that
$$
\|\Delta u_{\eps}\|_{L^\infty(\Omega,\RR^N )}+\|u_\vareps\|_{C^{1,\beta}(\Omega,\RR^N)} \leq C, \quad \forall~\eps >0.
$$
Moreover, $u_\vareps \rightarrow u_0$ in $C^{1,\beta'}(\Omega,\RR^N)$ as $\eps \rightarrow 0$ for any $\beta' \in (0,\beta)$ and $u_0 \in C^{1,\beta}(\Omega,\RR^N)$.
\end{theorem}

(See also Proposition \ref{pro:main1} for a localized version of Theorem \ref{thm:MDI}.)

To dispel confusion, the space $C^{1,\beta}(\Omega,\RR^N) = C^{1,\beta}(\bar\Omega,\RR^N)$ (with $\beta \in (0,1])$ in the above theorem consists of maps $u: \bar \Omega \rightarrow \RR^N$ such that 
\[
\|u\|_{C^{1,\beta}(\Omega,\RR^N)} := \|u\|_{C^1(\bar\Omega,\RR^N)} + [\nabla u]_{C^{0,\beta}(\Omega,\RR^{N\times M})} < \infty
\]
where
\begin{align*}
 \|u\|_{C^1(\bar\Omega,\RR^N)}
 	&= \sup_{\Omega}\, (|u| + |\nabla u|),\\
[\nabla u]_{C^{0,\beta}(\Omega,\RR^{N\times M})}
	& = \sup_{x, y \in \Omega, x \neq y} \frac{|\nabla u(x) - \nabla u(y)|}{|x - y|^\beta}.
\end{align*}

\begin{remark}
We highlight the fact that our theorem shows that every stationary harmonic map $u_0\in H^1(\Omega, \Ss^{N-1})\cap C^0(\bar \Omega, \Ss^{N-1})$ with $u_0|_{\partial\Omega} \in C^{1, \beta}(\partial \Omega, \Ss^{N-1})$ which is a limit of critical points of $\mcF_\eps$ belongs to $C^{1, \beta}(\Omega, \Ss^{N-1})$. We believe this is true without the assumption that it is a limit of critical points of $\mcF_\eps$, but we only managed to find some references in which it is shown that $u_0 \in C^{0,\gamma}(\Omega, \Ss^{N-1})$ for $\gamma \in (0,1)$. However, if the assumption of the continuity of $u_0$ up to $\partial\Omega$ is relaxed to only interior continuity (while still keeping $u_0|_{\partial\Omega} \in C^{1, \beta}(\partial \Omega, \Ss^{N-1})$), $u_0$ may not be regular up to $\partial\Omega$ -- see the example in \cite{Poon}.

\end{remark}

Some comments on the proof of Theorem \ref{thm:MDI} are in order. In view of the estimate $|u_\eps| \leq 1$ in Theorem \ref{thm:MPI}, we can adapt existing arguments to obtain a uniform $L^\infty$ estimate for $\nabla u_\eps$ -- see Sections \ref{Sec:MonF} and \ref{Sec:GradEst}. In those steps, condition \eqref{ass:Wnondeg} is not needed. By \eqref{eq:eps}, estimating $\Delta u_\eps$ is equivalent to estimating $ \eps^{-2}W'(1 - |u_\eps|^2)$ (uniformly in $\eps$). By \eqref{ass:Wnondeg}, this in turn is equivalent to estimating $ \eps^{-\frac{2(\alpha + 1)}{\alpha}}W(1 - |u_\eps|^2)$. To this end, we observe that $W(1 - |u_\eps|^2)$ satisfies a differential inequality of the form
\[
-\eps^2\Delta f + a f^{\frac{2\alpha}{\alpha + 1}} \leq C \quad \text{ in } \Omega
\]
for some positive constants $a > 0, C > 0$. When $\alpha = 1$, this is a linear inequality and can be estimated using Bessel functions. When $\alpha > 1$, this is a semi-linear inequality that resembles the Loewner-Nirenberg equation in conformal geometry \cite{LoewnerNirenberg}, and we borrow ideas from this line of work to obtain the desired bound. See Section \ref{Sec:LapEst}.

In our companion paper \cite{INSZ26-2nd}, we use Theorems \ref{thm:MPI} and \ref{thm:MDI} to prove the following existence and uniqueness result (that requires slightly stronger assumptions on the potential):

\begin{theorem}\label{Thm:Uniq}
Let $M \geq 1, N \geq 2$, $\Omega \subset \RR^M$ be a bounded domain with $C^2$ boundary $\partial\Omega$. Suppose $W \in C^2((-\infty, 1])$, $W(0) = 0, W'(0)=0$, $W(t) > 0$ for all $t \neq 0$, and there exist $\alpha \geq 1$ and $c_0 > 0$ so that
$$\lim_{t \rightarrow 0} \frac{|W''(t) - c_0 \alpha(\alpha+1) |t|^{\alpha-1}|}{|t|^{\alpha-1}}   = 0.$$
Let $u_0: \Omega \rightarrow \mathbb{S}^{N-1}$ be a harmonic map  such that 
$u_0\in C^{1, \beta_*}(\Omega,\mathbb{S}^{N-1})$, $0 < \beta_* \leq 1$ and the linearized harmonic map operator at $u_0$
\[
L_\parallel: \zeta\in \{H^1_0(\Omega, \RR^N)\, :\, \zeta\cdot u_0=0\} \mapsto -\Delta \zeta  - 2(\nabla u_0: \nabla \zeta) u_0 - |\nabla u_0|^2 \zeta\in H^{-1}(\Omega,\R^N)
\]
is injective.

 Then, for every $0 < \beta < \beta_*$, there exist $\delta_*>0$ and $\eps_*>0$, depending only on $M, N, \Omega, W, \|u_0\|_{C^{1,\beta}(\Omega,\RR^N)}$, $\beta$ and $\beta_*$ such that, for every $\eps\in (0, \eps_*)$, there exists a unique solution $u_\eps$ of the Ginzburg-Landau system \eqref{eq:eps} with boundary value $u_{bd}=u_0$ in the $C^{1, \beta}(\Omega)$-ball of radius $\delta_*$ around $u_0$.
 \end{theorem}

The rest of the paper is organized as follows. We start in Section \ref{Sec:MP} with the proof of Theorem \ref{thm:MPI}. We study the monotonicity formula in Section \ref{Sec:MonF} and establish gradient estimates in Section \ref{Sec:GradEst}. Theorem \ref{thm:MDI} is treated in Section \ref{Sec:LapEst}. We also include an Appendix in which we give some miscellaneous results on the weak and stationary harmonicity of the limit $u_0$ as well as some elliptic estimates used in the body of the paper.

%%%%%%%%%%%%%%%%%%%%%%%%%%
%%%%%%%%%%%%%%%%%%%%%%%%%%%%
%%%%%%%%%%%%%%%%%%%%%%%%%%%

\section{Maximum principle}\label{Sec:MP}
In this section, we prove Theorem \ref{thm:MPI}. Namely, we fix $\eps > 0$ and prove that $|u|\leq 1$ in $\Omega$ for any critical point $u\in H^1\cap L^\infty(\Omega, \R^N)$ of $\mcF_\eps$ such that $u\big|_{\partial \Omega}\in C^{0}(\partial \Omega, \Ss^{N-1})$. When $W$ satisfies $W'(t) \leq 0$ for $t < 0$, this is rather well-known. Under the general assumptions we make on $W$ in the present paper, this maximum principle for $|u|$ is new to our knowledge and its proof requires new ingredients.

As $W\in C^2((-\infty,1])$, a critical point $u\in H^1\cap L^\infty(\Omega, \R^N)$ of $\mcF_\eps$ such that $u\big|_{\partial \Omega}\in C^{0}(\partial \Omega, \Ss^{N-1})$ belongs to $C^2(\Omega) \cap C^0(\bar \Omega)$. We will prove that every solution $u\in C^2(\Omega) \cap C^0(\bar \Omega)$ of  the (more general) problem\footnote{A-priori, the boundary data could depend on $\eps$ in \eqref{Eq:EL}.}
\begin{equation}
\begin{cases}
-\Delta u   = \frac{1}{\vareps^2} W'(1-|u|^2)u &\text{ in } \Omega,\\
|u| \leq 1 & \text{ on } \partial \Omega
\end{cases}
\label{Eq:EL}
\end{equation}
satisfies $|u|\leq 1$ in $\Omega$. 

\begin{theorem}
\label{thm:maxG}
Let $W \in C^2((-\infty,1])$  satisfy \eqref{ass:Wsweak}. If $\vareps > 0$, $\Omega\subset \R^M$ is a bounded domain and
$u\in C^2(\Omega) \cap C^0(\bar \Omega)$ is a solution to \eqref{Eq:EL}, then $|u| \leq 1$ in $\Omega$.
\end{theorem}

The proof of Theorem \ref{thm:maxG} is done in several steps by analysing the modulus of $u$, i.e., 
\[
\rho = |u|
\]
(see Lemmas \ref{lem:eq-mod} and \ref{Lem:W'Sign}) 
combined with the behaviour of entire radial solutions to the Ginzburg-Landau problem \eqref{Eq:SGL} (see Lemma \ref{Lem:ShBar} below). We start by proving that $\rho$ is a sub-solution of \eqref{Eq:SGL}.

\begin{lemma}
\label{lem:eq-mod} Let $W\in C^2((-\infty,1])$, $\vareps > 0$ and $\Omega \subset \RR^M$ be a domain. If $u\in C^2(\Omega)$ satisfies $-\Delta u   = \frac{1}{\vareps^2} W'(1-|u|^2)u$ in $\Omega$, then the modulus $\rho=|u|$ belongs to
$C^{0,1}_{\rm loc}(\Omega)$ and satisfies in the sense of distribution
\begin{equation}
-\Delta \rho
	\leq \frac{1}{\vareps^2}  W'(1-\rho^2)  \rho \quad \text{ in } \Omega.
	\label{Eq:rho-w}
\end{equation}
\end{lemma}

\begin{proof}
We start by computing in the open set $\{\rho > 0\}$ where $\rho$ is $C^2$:
\begin{align*}
\Delta \rho
	&= \frac{1}{\rho } u \cdot \Delta u + \frac{|\nabla u|^2 - |\nabla \rho |^2}{\rho }
	= -\frac{1}{\vareps^2}  W'(1-\rho^2) \rho + \frac{|\nabla u|^2 - |\nabla \rho |^2}{\rho }.
\end{align*}
Noting that  $|\nabla \rho | \leq  |\nabla u|$ in $\{\rho > 0\}$, we deduce that
\begin{equation}
\Delta \rho
	\geq -\frac{1}{\vareps^2}  W'(1-\rho^2)  \rho  \quad \text{ in } \{\rho > 0\}.
	\label{Eq:rho}
\end{equation}
Finally, we extend this inequality in the whole domain $\Omega$. Indeed, as
\[
\big|\rho(x) - \rho(y)\big| \leq |u(x) - u(y)|,
\]
we have $\rho \in C^{0,1}_{\rm loc}(\Omega) $ and $|\nabla \rho| \leq |\nabla u|$ a.e. in $\Omega$. 
For $\delta > 0$, let $$\rho_\delta = (\delta^2 + |u|^2)^{1/2} \in C^2(\Omega) .$$ It is clear that $0 \leq \rho_\delta - \rho \leq \delta$ and so $\rho_\delta \rightarrow \rho$ uniformly in $\Omega$ as $\delta \rightarrow 0$. Noting that
\[
|\nabla \rho_\delta| \leq \frac{ |u| |\nabla u|}{(\delta^2 + |u|^2)^{1/2}} \leq |\nabla u|,
\]
we deduce that $\rho_\delta$ converges to $\rho$ weakly in $H^1_{\rm loc}(\Omega)$. Now, $\rho_\delta$ satisfies
\[
\Delta\rho_\delta
	= -\frac{1}{\vareps^2} W'(1 - \rho^2) \frac{\rho^2}{\rho_\delta} + \frac{|\nabla u|^2 - |\nabla \rho_\delta|^2}{\rho_\delta} \geq -\frac{1}{\vareps^2} W'(1 - \rho^2) \frac{\rho^2}{\rho_\delta} \text{ in } \Omega.
\]
Sending $\delta \rightarrow 0$, as $\Delta\rho_\delta\to \Delta\rho$ in the sense of distributions in $\Omega$, we obtain \eqref{Eq:rho-w}.
\end{proof}

\begin{lemma}\label{Lem:W'Sign}
Let $W\in C^2((-\infty,1])$, $\vareps > 0$ and $\Omega \subset \RR^M$ be a domain with non-empty boundary. If $u \in C^2(\Omega) \cap C^0(\bar \Omega)$ is a solution of \eqref{Eq:EL} such that $\rho = |u|$ attains an interior maximum at $x_0 \in \Omega$ with $\rho(x_0) > 1$, then
\[
W'(1 - \rho(x_0)^2) > 0.
\]
\end{lemma}

\begin{proof} Indeed, if $W'(1-\rho(x_0)^2) \leq 0$, then the constant function $\bar\rho = \rho(x_0)$ satisfies
\[
\Delta \bar\rho = 0 \leq -  \frac{1}{\vareps^2}  W'(1-\bar\rho^2)\bar\rho \text{ in } \Omega.
\]
Let $\omega$ be the connected component of $\{\rho > 0\}$ containing $x_0$. Recall that $\rho \in C^2(\omega) \cap C^0(\bar\omega)$. Since $\rho \leq \bar\rho$ and $\rho(x_0) = \bar\rho(x_0)$, we deduce from the above inequality, \eqref{Eq:rho} and the strong maximum principle that $\rho \equiv \bar\rho > 1$ on $\bar \omega$. This is absurd as $\rho  \leq 1$ on $\partial \Omega \neq \emptyset$ and $\rho = 0$ on $\partial\omega \cap \Omega$.
\end{proof}

Note that if the potential $W$ satisfies additionally that $W'(t) \leq 0$ for every $t < 0$, the above lemma proves that $|u| \leq 1$ in $\Omega$. The role of Theorem \ref{thm:maxG} is to prove $|u| \leq 1$ without such additional requirement on the sign of $W'$. Therefore, in the following, we focus on the situation where the potential $W$ has points $t< 0$ with $W'(t)>0$ by studying the behaviour of some entire radial solutions of the Ginzburg-Landau equation \eqref{Eq:SGL} in $\R^M$.

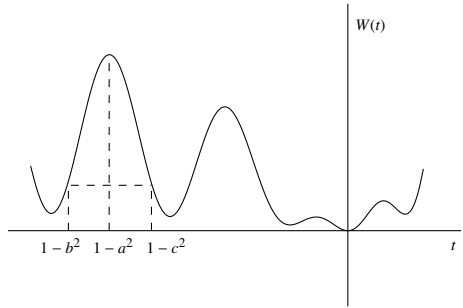
\begin{figure}[h]
\begin{center}
\begin{tikzpicture}
\draw (0,-1)--(0,3);
\draw (0.3,2.7) node  {\tiny $W(t)$};
\draw (-4.5,0)--(1.5,0);
\draw (1.4,-.2) node  {\tiny $t$};
\draw[smooth,samples=100,domain = -4.2:1] plot (\x,{1.5*\x*\x/(\x*\x-\x+1)*(1.2+cos(4*\x r))});
\draw[dashed]   (-3.7,0)--(-3.7,.6)--(-2.6,.6);
\draw (-3.8, -.2) node {\tiny $1 - b^2$};
\draw[dashed] (-2.60,0)--(-2.60,.6);
\draw (-2.4, -.2) node {\tiny $1 - c^2$};
\draw[dashed] (-3.16,0)--(-3.16,2.3);
\draw (-3.1, -.2) node {\tiny $1 - a^2$};
\end{tikzpicture}
\caption{The role of $a$, $b$, $c$ in Lemma \ref{Lem:ShBar}.}
\label{Fig1}
\end{center}
\end{figure}

\begin{lemma}\label{Lem:ShBar}
Suppose $W \in C^2((-\infty,1])$ satisfies \eqref{ass:Wsweak}. Let $1 < a < b \leq \infty$ be such that $W'(1 - s^2) > 0$ for $s \in (a,b)$, and let $c$ be the largest number in $[1,a)$ such that $W(1 - c^2) = W(1 - b^2)$ (see Figure \ref{Fig1}). For every $d \in [a,b)$ and $\vareps > 0$, there exists a unique solution $g = g_{\vareps,d} \in C^2([0,\infty))$ to the initial value problem
\[
\begin{cases}
g''(r) + \frac{M-1}{r} g'(r) = -\frac{1}{\vareps^2}  W'(1-g(r)^2) g(r) \quad \text{ in } (0,\infty),\\
g(0) = d, g'(0) = 0.
\end{cases}
\]
Moreover $g_{\vareps,d}$ satisfies $c < g_{\vareps,d} \leq d$ in $[0,\infty)$.
\end{lemma}

\begin{remark}\label{Rem:R1}
As a radial function on $\RR^M$, $g_{\vareps,d} = g_{\vareps,d}(|x|)$ satisfies \eqref{Eq:SGL} and $g_{\vareps,d} (0) = d$. As $c < g_{\vareps,d} \leq d$, we have by elliptic estimates that, for every $0< \vareps_1 < \vareps_2 < \infty$ and $R > 0$, the family $\{g_{\vareps,d}: \vareps \in [\vareps_1, \vareps_2], d \in [a,b)\}$ is bounded in $C^3([0,R])$. The uniqueness of $g_{\vareps,d}$ then implies that the map $(\vareps,d) \mapsto g_{\vareps,d}$ from $[\vareps_1,\vareps_2]  \times [a,b)$ into $C^2([0,R])$ is continuous.
\end{remark}

\begin{proof} We divide the proof in several steps:

\medskip

\nd {\it Step 1: We prove that $c <g \leq d$ in $[0,\infty)$ for any solution $g$.} For that, let $A(r) =  (g')^2 - \frac{1}{\vareps^2}   W(1 - g^2)$ for every $r\geq 0$. (See \cite{Modica-CPAM85}, where the function $A$ was used in the study of entire solutions to \eqref{Eq:SGL}.) By \eqref{Eq:SGL}, we compute
\begin{align*}
A'(r) 
	&= 2   g' g''  + \frac{2}{\vareps^2}   W'(1 - g^2) g g'
	= - \frac{2(M-1)}{r}  (g')^2  \leq 0 \quad \text{ in } (0,\infty).
\end{align*}
As $g(0)=d$ and $g'(0)=0$, it follows that $A(r) \leq A(0) = - \frac{1}{\vareps^2}   W(1 - d^2)$ in $[0,\infty)$. This implies that 
\begin{equation}
W(1 - g^2) \geq W(1 - d^2) \quad \text{ in }[0,\infty). 
\label{Eq:Wgd}
\end{equation}
By the hypothesis that $W'(1 - s^2) > 0$ for $s \in (d,b) \subset (a,b)$, we have
\begin{equation}
W(1 - b^2) < W(1 - s^2) < W(1 - d^2)
	\quad \text{ for } s \in (d,b).
\label{Eq:Wdb}
\end{equation}
It follows from \eqref{Eq:Wgd} and \eqref{Eq:Wdb} that $g$ cannot attain any value in $(d,b)$. By the continuity of $g$ and the fact that $g(0) = d$, we must have $g \leq d$ in $(0,\infty)$. 
As $W(1 - g^2) \geq W(1 - d^2)> W(1 - b^2)=W(1-c^2)$ and $g(0)=d>c$, by continuity of $g$, we conclude that $g > c$ in $[0,\infty)$.
\medskip

\noindent {\it Step 2: We prove the existence of a solution $g$.}
For $\delta > 0$, let $g_\delta$ be the solution to 
\[
\begin{cases}
g_\delta''(r) + \frac{M-1}{r} g_\delta'(r) = -\frac{1}{\vareps^2}  W'(1-g_\delta^2) g_\delta \quad \text{ for } r > \delta,\\
g_\delta(\delta) = d, g_\delta'(\delta) = 0.
\end{cases}
\]
Standard ODE theories imply that $g_\delta$ exists for some maximal interval $[\delta, R)$ with $R \leq \infty$. Arguing as in Step 1, with $A$ replaced by $A_\delta(r) =  (g_\delta')^2 - \frac{1}{\vareps^2}   W(1 - g_\delta^2)$ and noting that $g_\delta(\delta) = d$ and $A_\delta(\delta) = - \frac{1}{\eps^2}W(1 - d^2)$, we have that 
\[
c < g_\delta \leq d \quad \text{ in } [\delta,R).
\]
Moreover, by the ODE, we also have $|(r^{M-1}g_\delta')'|\leq Cr^{M-1}$ for every $r\in [\delta, R)$ and by integration, as $g_\delta'(\delta)=0$, we get $|g_\delta'(r)|\leq \tilde C r$ for every $r\in [\delta, R)$. (The constants $C$ and $\tilde C$ may depend on $\vareps$, which is fixed in this proof.) Therefore, $R = \infty$. This also implies that the family $(g_\delta)_{\delta \downarrow 0}$ is bounded in $C^3(K)$ for any compact set $K \subset (0,\infty)$. Thus, along a sequence $\delta_m \rightarrow 0$, $g_{\delta_m}$ converges in $C^2_{\rm loc}((0,\infty))$ to some $g$ satisfying $c \leq g \leq d$ and
\[
g''(r) + \frac{M-1}{r} g'(r) = -\frac{1}{\vareps^2}  W'(1-g^2) g \quad \text{ in } (0,\infty).
\]
It remains to show that $g$ extends to a function in $C^2([0,\infty))$ with $g(0) = d$ and $g'(0) = 0$. Indeed, we have for every $r > \delta_m$
\begin{align*}
g_{\delta_m}'(r) 
	&= -\frac{1}{\vareps^2 r^{M-1}} \int_{\delta_m}^r s^{M-1} W'(1-g_{\delta_m}(s)^2) g_{\delta_m}(s)\, ds,\\
g_{\delta_m}(r) 
	&= d + \int_{\delta_m}^r g_{\delta_m}'(s)\, ds.
\end{align*}
Thus, if we fix some $Q \geq \max_{t \in [c,d]} |W'(1 - t^2)| t$, then
\begin{align*}
|g_{\delta_m}'(r)| &\leq \frac{Q}{M \vareps^2 }\, r , \quad
|g_{\delta_m} - d| \leq \frac{Q}{2M \vareps^2 }\, r^2 \quad \text{ for } r \in [\delta_m,\infty).
\end{align*}
Passing to the limit $\delta_m\to 0$, this implies that 
\begin{align*}
|g'(r)| &\leq \frac{Q}{M \vareps^2 }\, r,\quad
|g(r) - d| \leq \frac{Q}{2M \vareps^2 }\, r^2 \quad \text{ for } r \in (0,\infty).
\end{align*}
In particular, $g$ is differentiable up to $r = 0$, $g(0) = d$ and $g'(0) = 0$. Recall that, as a radial function on $\RR^M$, $g \in C^1(\RR^M)$ and $g$ satisfies
\[
\Delta g = -\frac{1}{\vareps^2}  W'(1-g^2) g \quad \text{ in } \RR^M \setminus \{0\}.
\]
Standard elliptic theories implies that $0$ is a removable singularity for $g$, $g \in C^2(\RR^M)$, and $g$ satisfies \eqref{Eq:SGL}. In particular, as a function of $r$, $g \in C^2([0,\infty))$.

\medskip

\noindent {\it Step 3: We prove uniqueness of $g$.} We only need to show the uniqueness of $g$ in some interval $[0,\mu]$ with $\mu > 0$. Once this is done, standard uniqueness theory for ODEs gives the uniqueness of $g$ in $[0,\infty)$.
Indeed, assume that $g_1$ and $g_2$ are two solutions in $[0,\mu]$ to our initial value problem. Then we have for $r \in [0,\mu]$:
\[
g_i  (r) = d -\frac{1}{\vareps^2} \int_0^r \frac{1}{s^{M-1}} \int_0^s \xi^{M-1} W'(1 - g_i(\xi)^2) g_i(\xi)\,d\xi\,ds, \quad i = 1,2.
\]
By Step 1, we have $c < g_i \leq d$; as $W \in C^2((-\infty,1])$, we have after possibly slightly enlarging the constant $Q$ introduced in Step 2 that
\[
|g_1(r) - g_2(r)| \leq \frac{Q r^2}{2M\vareps^2} \sup_{\xi \in [0,\mu]} | g_1(\xi) - g_2(\xi) | , \quad \forall r\in [0, \mu].
\]
Clearly, when $\frac{Q}{2M\vareps^2}\mu^2 < 1$, this is possible if and only if $g_1 = g_2$ in $[0,\mu]$. The proof is complete.
\end{proof}

Now, we prove the desired maximum principle for $|u|$ stated in Theorem \ref{thm:maxG}:

\begin{proof}[Proof of Theorem \ref{thm:maxG}]
Suppose by contradiction that $\rho = |u|$ satisfies
\[
m = \max_{\bar \Omega} \rho > 1.
\]
(Note that $m$ is a finite number as $\Omega$ is bounded.) Then $\rho \in [0,m]$ in $\bar\Omega$. Since $\rho \leq 1$ on $\partial \Omega$, $m$ is attained in $\Omega$. By Lemma \ref{Lem:W'Sign}, $W'(1 - m^2) > 0$. Define
\begin{align*}
a &= \inf\Big\{1\leq t < m: W'(1 - s^2) > 0 \text{ for every } s\in (t, m) \Big\} \in [1, m),\\
b &= \sup\Big\{t > m: W'(1 - s^2) > 0 \text{ for every } s\in (m, t)\Big\} \in (m, \infty].
\end{align*}
Clearly, $a$ and $b$ are well-defined by continuity of $W'$. Since $W'(1-m^2) > 0$ and $0$ is a global minimum of $W$ (by \eqref{ass:Wsweak}), $W$ has a local maximum in $(1-m^2,0)$. Hence $a > 1$. Also, by adjusting $W$ in $(-\infty, -C)$ for some large $C>0$ so that $W' < 0$ in $(-\infty, -C)$ while keeping $W$ the same in $[1 - m^2, 1]$ so that the equation \eqref{Eq:EL} for $u$ remains unchanged, we may assume that $b < \infty$. We then have $W'(1 - t^2) > 0$ for $a < t < b$, $W'(1 - a^2) = W'(1 - b^2) = 0$. Let $c$ be the  largest number in $[1,a)$ such that $W(1 - c^2) = W(1 - b^2)$. For $d \in [a,b)$, let $g_{\vareps,d}$ be given by Lemma \ref{Lem:ShBar} so that, as a radial function in $\RR^M$, it satisfies $c < g_{\vareps,d} \leq d$ and
\[
\begin{cases}
\Delta g_{\vareps,d} = -\frac{1}{\vareps^2}  W'(1-g_{\vareps,d}^2) g_{\vareps,d} \quad \text{ in } \RR^M,\\
g_{\vareps,d}(0) = d.
\end{cases}
\]
By uniqueness, we have $g_{\vareps, a} \equiv a$. We claim that
\[
\lim_{d \nearrow b} g_{\vareps,d} = b,
\]
where the limit is uniform on compact subset of $[0,\infty)$. Indeed, by Remark \ref{Rem:R1}, for any sequence $d_j \nearrow b$, there is a subsequence $d_{j_k}$ such that $g_{\vareps, d_{j_k}}$ converges in $C^2_{\rm loc}([0,\infty))$ to a function $g_*$ satisfying $c \leq g_* \leq b$ and (after viewing $g_*$ as a radial function in $\RR^M$),
\[
\begin{cases}
\Delta g_* = -\frac{1}{\vareps^2}  W'(1-g_*^2) g_* \quad \text{ in } \RR^M,\\
g_*(0) = b.
\end{cases}
\]
As the constant function $\bar g = b$ satisfies $\Delta \bar g = -\frac{1}{\vareps^2}  W'(1-\bar g^2) \bar g$ in $\RR^M$, $\bar g \geq g_*$ and $\bar g(0) = g_*(0)$, we deduce from the strong maximum principle that $g_* \equiv \bar g = b$. The claim follows.

By the claim and the fact that $m = \max \rho < b$ together with the boundedness of $\Omega$, there exists some small $\mu > 0$ such that $\rho < g_{\vareps,d}$ in $\Omega$ for $b - \mu < d < b$ (after viewing $g_{\vareps,d}$ as a radial function in $\RR^M$). Define
\[
\underline{d} = \inf\Big\{a \leq d < b:\rho < g_{\vareps,d} \text{ in } \Omega\Big\}.
\]
By Remark \ref{Rem:R1}, we have $\rho \leq g_{\vareps,\underline{d}}$ in $\Omega$. First, $\underline{d}>a$; because  otherwise, $\underline{d}=a$ and $\rho \leq g_{\vareps,a} = a < m$ in $\Omega$, which is a contradiction to the definition of $m$. Second, by the definition of $\underline{d}$, there exists $x_1 \in  \bar \Omega$ such that $\rho(x_1) = g_{\vareps,\underline{d}}(|x_1|) > c \geq 1$ (otherwise $\underline{d}$ could be decreased contradicting its definition) where we used the lower bound in Lemma \ref{Lem:ShBar} for $ g_{\vareps,\underline{d}}$. As $\rho  \leq 1$ on $\partial \Omega$, we get $x_1\in \Omega$. Let $\omega$ be the connected component of $\{\rho > 0\}$ containing $x_1$. Noting that $\rho \in C^2(\omega) \cap C^0(\bar\omega)$, in view of \eqref{Eq:rho} and the equation for $g_{\vareps, \underline{d}}$, we can apply the strong maximum principle to deduce that $\rho \equiv g_{\vareps,\underline{d}}$ in $\omega$. We then reach a contradiction with the fact that $\rho  \leq 1 \leq c < g_{\vareps,\underline{d}}$ on $\partial \Omega \neq \emptyset$ and $\rho = 0 < c < g_{\vareps,\underline{d}}$ on $\partial \omega \cap \Omega$. The proof is completed.
\end{proof}

\begin{proof}[Proof of Theorem \ref{thm:MPI}] 
As $W'\in C^1((-\infty,1])$ and $u\in L^\infty(\Omega)$, we have $\Delta u\in L^\infty(\Omega)$ and by De Giorgi-Nash-Moser's theorem (e.g. \cite[Theorem 8.24]{GT}) that $u \in C^{0,\gamma}_{\rm loc}(\Omega)$ for some $\gamma \in (0,1)$. It follows that $-\Delta u = \frac{1}{\vareps^2} W'(1 - |u|^2) u \in  C^{0,\gamma}_{\rm loc}(\Omega)$. By interior Schauder's estimate, we have $u \in C^{2,\gamma}_{\rm loc}(\Omega)$ (see \cite[Theorem 4.13]{GT}). Since $u\big|_{\partial \Omega}\in C^0(\partial \Omega)$, $\Delta u\in L^\infty(\Omega)$ and $\Omega$ satisfies an exterior cone condition at each point on $\partial\Omega$, we then have  that $u \in C^0(\bar \Omega)$ (see \cite[Theorem 8.30]{GT}). We can apply Theorem \ref{thm:maxG} to obtain the desired conclusion.
\end{proof}

Note that in Theorem \ref{thm:MPI}, the condition \eqref{ass:Wsweak} cannot be dropped:
 \begin{proposition}\label{Prop:WsweakCE}
If $W \geq 0$ in $(0,1)$ and $W$ is negative somewhere in $(-\infty,0)$, then for any given boundary data $u_{bd} \in H^{1/2}(\partial\Omega,\Ss^{N-1})$, when $\eps$ is sufficiently small, any minimizer $u_\eps\in H^1 \cap L^\infty(\Omega,\RR^N)$  of $\mcF_\eps$ with the boundary data $u_{bd}$, if it exists, must satisfy $|u_\eps| > 1$ somewhere in $\Omega$.
\end{proposition}

\begin{proof}
Since $\mcF_\eps[u] \geq 0$ for all $u \in H^1(\Omega, \RR^N)$ with $|u| \leq 1$, it suffices to exhibit some $v \in H^1(\Omega,\RR^N)$ with $v = u_{bd}$ on $\partial\Omega$ such that $\mcF_\eps[v] < 0$ for all sufficiently small $\eps > 0$. For some small $\delta > 0$ to be chosen later, let $\Omega_\delta= \{x \in \Omega\, : \,  \mathrm{dist}(x,\partial\Omega) > \delta\}$ and select a cut-off function $\chi \in C_c^{0,1}(\Omega)$ such that $0 \leq \chi \leq 1$ in $\Omega$ and $\chi \equiv 1$ in $\Omega_\delta$. Fix a vector $e \in \RR^N$ such that $W(1 - |e|^2) < 0$ (such $e$ exists since $W$ is negative somewhere in $(-\infty, 0)$). Define $w$ as the harmonic extension of $u_{bd}$ in $\Omega$, in particular, $w\in H^1 \cap L^\infty(\Omega,\RR^N)$; then set
\[
v = (1 - \chi)w  + \chi e \in H^1 \cap L^\infty(\Omega,\RR^N).
\]
Note that $\|v\|_{L^\infty(\Omega)} \leq \|w\|_{L^\infty(\Omega)} + |e| =: m$. We have
\begin{align*}
\int_\Omega W(1 - |v|^2)\,dx 
	&= \int_{\Omega \setminus \Omega_\delta}W(1 - |v|^2)\,dx + \int_{\Omega_\delta}W(1 - |v|^2)\,dx \\
	&\leq |\Omega \setminus \Omega_\delta| \max_{[1 - m^2,1]} W 
		+ |\Omega_\delta| W(1 - |e|^2) .
\end{align*}
As $W(1 - |e|^2) < 0$ and $ |\Omega \setminus \Omega_\delta| \to 0$ as $\delta\to 0$, we may choose a sufficiently small $\delta > 0$ so that 
\[
\int_\Omega W(1 - |v|^2)\,dx < 0.
\]
It follows that
\[
\mcF_\eps(v) = \frac{1}{2}\int_\Omega |\nabla v|^2 + \frac{1}{2\eps^2}  \int_\Omega W(1 - |v|^2)\,dx < 0 \text{ for all sufficiently small } \eps > 0.
\]
This entails that $\mcF_\eps(u_\eps)<0$ so $|u_\eps|$ achieves values $>1$.
\end{proof}

\begin{remark} If in addition to the assumptions of the above proposition, there exists a constant $C>0$ such that $W'(t) \leq 0$ for $t \leq - C$, then there exists a minimizer $u_\eps\in H^1 \cap L^\infty(\Omega,\RR^N)$ of $\mcF_\eps$ for any given boundary data $u_{bd} \in H^{1/2}(\partial\Omega,\Ss^{N-1})$. Indeed, this is due to a cutting-off procedure: for every map $u\in  H^1 (\Omega,\RR^N)$, one assigns the map $$v=\begin{cases} u & \textrm{if } |u|^2\leq C+1\\ \sqrt{C+1}\frac{u}{|u|} & \textrm{if } |u|^2\geq C+1\end{cases}$$ yielding  $\mcF_\eps[v]\leq \mcF_\eps[u]$, $v \in H^1 \cap L^\infty(\Omega,\RR^N)$ and $v=u_{bd}$ on $\partial \Omega$; then, by the direct method, one gets the existence of a bounded minimizer.
\end{remark}

We end this section with a remark of higher regularity for $u$ if the boundary $u_{bd}$ has better regularity. 

\begin{remark}
\label{rem:regu}
If $W \in C^2$, $\partial\Omega$ is $C^{1,\beta}$ with $\beta \in (0,1)$ and $u$ is a critical point $\mcF_\vareps$ in $H^1\cap L^\infty(\Omega, \R^N)$ subject to the boundary data $u_{bd}\in 
C^{1,\beta}(\partial \Omega, \Ss^{N-1})$, then, by elliptic estimates, we have $u \in C^2(\Omega) \cap C^{1,\beta}( \Omega)$  (see \cite[Corollary 8.36]{GT}; see also \cite[Theorems 4.15, 4.16]{GT}). Moreover, $|\nabla^2 u (x)| \leq C_{u,\eps} d(x,\partial\Omega)^{\beta - 1}$ in $\Omega$ (to see this, apply a rescaled version of \cite[Theorem 4.6]{GT} to the function $y \mapsto u(y) - u(x) - \nabla u(x) \cdot (y - x)$ in the ball $B_{d(x ,\partial\Omega)/2}(x )$). In particular, $u  \in W^{2,1}(\Omega)$.
\end{remark}

%%%%%%%%%%%%%%%

%%%%%%%%%%%%%
\section{Monotonicity formula}\label{Sec:MonF}

Let $B_r(x_0)$ be the ball centered at $x_0$ of radius $r$. The aim of this section is to prove, for a critical point $u$ of the energy $\mcF_\eps$ defined by \eqref{def:mcFeps}, a monotonicity formula for the following quantity corresponding to the concentration of the energy on sets of codimension $2$: for $x_0 \in \bar\Omega$ and $r > 0$, let 
\[
\varphi_{x_0}(r) = \frac{1}{r^{M-2}} \int_{\Omega \cap B_r(x_0)} e_\eps(u)\,dx,
\]
where $e_\eps(u)$ is the energy density in \eqref{def:mcFeps}: 
\be\label{def:edensity}
e_\eps(u):=\frac{1}{2}|\nabla u|^2+\frac{1}{2\eps^2}W(1-|u|^2).
\ee
When $x_0$ is the origin, we will write $\varphi$ and $B_r$ in place of $\varphi_0$ and $B_r(0)$.

We establish:

\begin{proposition}\label{prop:monot}
Let $W \in C^2((-\infty,1])$ satisfy \eqref{ass:Wsweak}, $\Omega\subset \R^M$ be a domain, $\eps > 0$ and $u\in H^1\cap L^\infty(\Omega, \R^N)$ be a critical point of $\mcF_\eps$.

\begin{enumerate}
\item For every $x_0 \in \Omega$, it holds that $\varphi_{x_0}$ is differentiable in the interval $(0, d(x_0,\partial\Omega))$  and 
 $
\varphi_{x_0}' \geq 0.$ In particular, 
$\varphi_{x_0}(r)\le \varphi_{x_0}(s)$ for all $r\le s<d(x_0,\partial\Omega)$.

\item Suppose further that $\Omega$ is a $C^2$ bounded domain and 
$u_{bd}\in C^{1, \beta}(\partial \Omega, \Ss^{N-1})$, $\beta \in (0,1)$, and $u = u_{bd}$ on $\partial\Omega$. Then there exist $r_0 > 0$ and $C > 1$ depending only on $\partial\Omega$ and $u_{bd}$ (and hence independent of $\eps$) such that, for every $x_0 \in \partial\Omega$, $\varphi_{x_0}$ is differentiable in $(0 , r_0 )$ and 
\be
\label{ineg22}
\varphi_{x_0}'(r) \geq - Cr - C\varphi_{x_0}(2r) - C\varphi_{x_0}(2r)^{1/2}, \quad \forall \rho\in (0 ,r_0/2).
\ee
As a consequence, for any $\delta \in (0,1)$, there exists $r_* > 0$ depending only on $\partial\Omega$, $u_{bd}$ and $\delta$ such that if $x_0\in \partial \Omega$, $0 < r\leq r_*$ and $\varphi_{x_0} \leq \delta$ in $[r/2,r]$, then $\varphi_{x_0} \leq 2\delta$ in $(0,r]$.
\end{enumerate}

\end{proposition}

As in works of many other authors, one important application of the above monotonicity formula is to deduce smallness of the energy at all scales from smallness of the energy at one fixed scale. 

The proof of Proposition \ref{prop:monot} uses available ideas in the  literature with some small tweaks. We however note that, in the case $x_0$ lies on the boundary, estimate \eqref{ineg22} is slightly different from some other similar estimates in the literature e.g. in \cite[Proposition 2.1]{CLR18}, \cite[Lemma 9]{MZ10} where the related estimates were of the form
\[
\varphi_{x_0}'(r) \geq - C - C\varphi_{x_0}(2r).
\]
Compare also \cite[Lemma II.5]{LinRiv2} where the boundary map was assumed to be in the stronger class $C^2$. 

To prove Proposition \ref{prop:monot}, we start with the following Pohozaev type identity for every smooth enough map $u$ (not necessarily a solution of \eqref{eq:eps}):

\begin{lemma}\label{Lem:VFM} 
Let $X \in C^1(\bar\Omega,\RR^M)$, $W\in C^2((-\infty, 1])$ and 
$u\in W^{2,1}(\Omega)\cap  C^{0,1}(\Omega)$ in a domain $\Omega\subset \R^M$. For any bounded subdomain $\omega \subset \Omega$ with Lipschitz boundary $\partial\omega$, we have
\begin{align*}
&\int_{\omega} \Big[-\Delta u - \frac{1}{\vareps^2} W'(1 - |u|^2)u\Big] \cdot [(X \cdot \nabla) u]\, dx\\
	&\quad=  \int_{\partial\omega}  \Big[-(\nu \cdot \nabla) u  \cdot (X \cdot \nabla) u + e_\eps(u) \nu \cdot X \Big]\,d\sigma
	+ \int_{\omega} \Big\{ \sum_{i,j=1}^M  \nabla_i X_j \nabla_i u\cdot \nabla_j u - e_\eps(u) \nabla \cdot X\Big\}\,dx,
\end{align*}
where $\nu$ is the unit outer normal vector field to $\partial \omega$.
\end{lemma}

\begin{proof} Note that
\begin{align*}
\Delta u\cdot [(X \cdot \nabla) u]
	&= \sum_{i,j,k} \Big\{\nabla_i \big(\nabla_i u_k X_j \nabla_j u_k\big) - \nabla_i u_k \nabla_i X_j \nabla_j u_k - \nabla_i u_k X_j \nabla_i \nabla_j u_k\Big\}\\
	&= \sum_{i,j,k} \Big\{\nabla_i \big(\nabla_i u_k X_j \nabla_j u_k\big) - \frac{1}{2} \nabla_j\big((\nabla_i u_k)^2 X_j\big)
		- \nabla_i u_k \nabla_i X_j \nabla_j u_k + \frac{1}{2} (\nabla_i u_k)^2 \nabla_j X_j \Big\},
\end{align*}
and
\[
- W'(1 - |u|^2)u \cdot [(X \cdot \nabla) u] = \frac{1}{2} X \cdot \nabla [W(1 - |u|^2)].
\]
Therefore
\begin{align*}
&\int_{\omega} \Big[-\Delta u - \frac{1}{\vareps^2} W'(1 - |u|^2)u\Big] \cdot [(X \cdot \nabla) u]\\
	&=   \int_{\partial\omega }  \Big[-(\nu \cdot \nabla) u  \cdot (X \cdot \nabla) u + \frac{1}{2}|\nabla u|^2 \nu \cdot X + \frac{1}{2\vareps^2} W(1 - |u|^2) \nu \cdot X\Big]\,d\sigma\\
		&\qquad + \int_{\omega} \Big\{ \sum_{i,j,k} \nabla_i u_k \nabla_i X_j \nabla_j u_k - \frac{1}{2} |\nabla u|^2 \nabla \cdot X - \frac{1}{2\vareps^2} W(1 - |u|^2) \nabla \cdot X\Big\}\,dx.
\end{align*}
Recalling the expression for $e_\eps(u)$, we obtain the lemma.
\end{proof}

\begin{corollary}[Boundary Rellich-type estimate]\label{Cor:BRE} 
Let $W \in C^2((-\infty,1])$ satisfy \eqref{ass:Wsweak}, $\Omega\subset \R^M$ be a bounded $C^2$ domain and 
$u_{bd}\in C^{1, \beta}(\partial \Omega, \Ss^{N-1})$, $\beta \in (0,1)$. Then there exist $r_0 > 0$ and $C > 1$ depending only on $\Omega$ such that for every $\eps>0$ and every critical point $u\in H^1\cap L^\infty(\Omega, \R^N)$ of $\mcF_\eps$ with the boundary condition $u_{bd}$, it holds:
\[
\int_{\partial\Omega \cap B_{r/2}(x_0)} | (\nu \cdot \nabla) u  |^2 \,d\sigma
	\leq C\int_{\partial\Omega \cap B_r(x_0)} |\nabla_{\partial\Omega} u|^2\,d\sigma
		+ \frac{C}{r}\int_{\Omega \cap B_r(x_0)}   e_\eps(u) \,dx
\]
for all $x_0 \in \partial\Omega$ and $0 < r < r_0$. Here, $\nu$ is the unit outer normal vector field to $\partial \Omega$ and $\nabla_{\partial \Omega}=\nabla -(\nu\otimes \nu) \nabla$ is the gradient in the tangential directions at $\partial \Omega$.

\end{corollary}

\begin{proof}
We fix $r_1 > 0$ sufficiently small (depending only on $\partial\Omega$) such that, for every $x \in \partial\Omega$, $\partial\Omega \cap B_{r_1}(x)$ is a graph of a $C^2$ function over a domain $U_x$ on the tangent plane of $\partial\Omega$ at $x$. Fix $0 < r_0 < r_1$ such that, for every $x \in \partial\Omega$ and $0 < r < r_0$, $\Omega \cap B_{r}(x)$ is a Lipschitz domain and its projection onto the tangent plane of $\partial\Omega$ at $x$ is a subset of $U_x$. We can assume without loss of generality that $x_0 = 0$. By a rotation of coordinate axes, we write $\partial \Omega \cap B_{r_0} = \{(x',G(x')) \in \RR^M: x' \in U = U_{0} \subset \RR^{M-1}\}$ and $\Omega \cap B_{r_0} = \{(x',x_M) \in B_{r_0}: x_M > G(x')\}$ where $G:U\to \R$ is a $C^2$ function satisfying $G(0) = 0, \nabla G(0) = 0$. Define the $C^1$ vector field $Y$ in $U \times \RR$ by $Y(x',x_M) = (\nabla G(x'),-1)$. Note that, along $\partial \Omega \cap B_{r_0}$, the unit normal $\nu$   is given by $\nu = \frac{Y}{|Y|}$.\footnote{ \label{Foot:nx}For our later use, we note that $\nu \cdot x = (1 + |\nabla G(x')|^2)^{-1/2} (\nabla G(x') \cdot x' - G(x')) \geq - C|x'|^2 \geq - C|x|^2$ along $\partial \Omega \cap B_{r_0}$ for some constant $C$ depending only on the $C^2$ norm of $G$.}

Let $\eta \in C_c^\infty([0,\infty))$   such that $0\leq \eta \leq 1$ in $[0,\infty)$, $\eta \equiv 1$ in $[0,1/2)$ and $\eta \equiv 0$ in $[1,\infty)$ and let $\eta_r(x) = \eta(|x|/r)$ for every $0<r<r_0$. In particular, $\supp \eta_r\subset \bar B_r$. Recall the regularity of our critical point $u\in W^{2,1}(\Omega)\cap  C^{1,\beta}(\Omega)$ (see Remark \ref{rem:regu}). Applying Lemma \ref{Lem:VFM} with the $C^1$ vector field $X = \eta_r Y$ and the subdomain $\Omega \cap B_r$ (noting that $W(1- |u|^2) = W(0) = 0$ on $\partial \Omega$), we have 
\begin{align*}
0 
&=  \frac{1}{2} \int_{\partial\Omega \cap B_r} \eta_r |Y| \Big[|\nabla_{\partial\Omega} u|^2- | (\nu \cdot \nabla) u  |^2 \Big]\,d\sigma\\
		&\qquad + \int_{\Omega \cap B_r} \Big\{ \sum_{i,j,k} \nabla_i u_k \nabla_i (\eta_r Y_j) \nabla_j u_k - e_\eps(u) \nabla \cdot (\eta_r Y)\Big\}\,dx.
\end{align*}
Since $e_\vareps(u) \geq \frac{1}{2} |\nabla u|^2 \geq 0$ (as $W \geq 0$) and $ |Y| \geq 1$ in $\Omega \cap B_{r_0}$, it follows that
\begin{align*}
\int_{\partial\Omega \cap B_{r/2}} | (\nu \cdot \nabla) u  |^2 \,d\sigma
	&\leq \int_{\partial\Omega \cap B_r} \eta_r |Y| | (\nu \cdot \nabla) u  |^2 \,d\sigma\\
	&\leq C\int_{\partial\Omega \cap B_r} |\nabla_{\partial\Omega} u|^2\,d\sigma
		+ \frac{C}{r}\int_{\Omega \cap B_r}   e_\eps(u) \,dx.
\end{align*}
The proof is complete.  
\end{proof}

Now, we prove the desired monotonicity formula:

\begin{proof}[Proof of Proposition \ref{prop:monot}] 
Let us first prove point 1. of the proposition. Let $x_0 \in \Omega$ and $0<r < d(x_0,\partial\Omega)$. Up to a translation, we may assume that $x_0=0$. By the definition of $\varphi(r):=\varphi_{x_0}(r)$ and the fact that $u\in C^{2}(\Omega)$ (by elliptic regularity), $\varphi$ is differentiable for $0<r < d(x_0,\partial\Omega)$. Applying Lemma \ref{Lem:VFM} with the identity vector field $X(x) = x$,  on $B_r \subset \Omega$, we have
\begin{align*}
0 
	&= r \int_{\partial B_r}  \big[-|\partial_r u|^2 + e_\eps(u)  \big]\,d\sigma
	 - \int_{B_r} \Big\{(M-2)e_\eps(u)  + \frac{1}{\vareps^2} W(1 - |u|^2)\Big\}\,dx.
\end{align*}
Therefore, for $0<r < d(x_0,\partial\Omega)$,
\begin{align*}
\varphi'(r) 
	&=   \frac{1}{r^{M-2}}\int_{\partial B_r} e_\eps(u)\,d\sigma-\frac{M-2}{r^{M-1}}\int_{B_r} e_\eps(u)\,dx\\
	&= \frac{1}{r^{M-2}} \int_{\partial B_r}  |\partial_r u|^2 \,d\sigma 
		+  \frac{1}{r^{M-1}}\int_{B_r} \frac{1}{\vareps^2} W(1 - |u|^2)\,dx \geq 0,
\end{align*}
where we used that $W\geq 0$. This concludes the proof of point 1.

Next, let us first prove point 2. of the proposition. Recall the regularity of our critical points $u\in W^{2,1}(\Omega)\cap  C^{1,\beta}(\Omega)$ (see Remark \ref{rem:regu}). Let $r_0 > 0$ be as in Corollary \ref{Cor:BRE}. Fix some $x_0 \in \partial\Omega$, which we assume without loss of generality to be the origin, and $0<r < r_0/2$. As $u\in C^{1,\beta}(\Omega)$ and recalling the way $r_0 < r_1$ were set up in the proof of Corollary \ref{Cor:BRE}, we have that $\varphi$ is differentiable for $0<r < r_0$. Applying Lemma \ref{Lem:VFM} with the identity vector field $X(x) = x$ and the subdomain $B_r \cap \Omega$ where $\nu$ is the unit outer normal field to $\partial(B_r \cap \Omega)$  (noting that $W(1- |u|^2) = W(0) = 0$ on $\partial \Omega$) , we deduce
\begin{align*}
\varphi'(r) 
	&= \frac{1}{r^{M-2}}\int_{\Omega \cap \partial B_r} e_\eps(u)\,d\sigma-\frac{M-2}{r^{M-1}}\int_{\Omega \cap B_r} e_\eps(u)\,dx \\
	&= \frac{1}{r^{M-2}} \int_{\Omega \cap  \partial B_r}  |\partial_r u|^2 \,d\sigma 
		+  \frac{1}{r^{M-1}}\int_{\Omega \cap B_r} \frac{1}{\vareps^2} W(1 - |u|^2)\,dx\\
		&\qquad +\frac{1}{r^{M-1}} \int_{\partial \Omega \cap B_r}  \Big[(\nu \cdot \nabla) u  \cdot (x \cdot \nabla) u - \frac{1}{2} |\nabla u|^2 \nu \cdot x \Big]\,d\sigma.
\end{align*}
The first two integrals are non-negative (recall that $W \geq 0$) and will be discarded to obtain a lower bound for $\varphi'(r)$. To treat the last integral, we have for $x\in \partial \Omega \cap B_r$:
\[
|(x \cdot \nabla)u - x \cdot \nu (\nu \cdot \nabla) u| = \left|\big([x - (x \cdot \nu)\nu] \cdot \nabla\big) u\right| \leq |x| |\nabla_{\partial\Omega} u_{bd}|,
\]
which after multiplying by $|(\nu \cdot \nabla) u|$ implies:
\[
(\nu \cdot \nabla) u  \cdot (x \cdot \nabla) u \geq |(\nu \cdot \nabla) u|^2 \nu \cdot x - |x| |(\nu \cdot \nabla) u|  |\nabla_{\partial\Omega} u_{bd}|.
\]

As $|\nabla u|^2=|\nabla_{\partial\Omega} u_{bd}|^2+|(\nu \cdot \nabla) u|^2$ on $\partial \Omega \cap B_r$, we obtain
\[
(\nu \cdot \nabla) u  \cdot (x \cdot \nabla) u - \frac{1}{2} |\nabla u|^2 \nu \cdot x
	\geq \frac{1}{2} |(\nu \cdot \nabla) u|^2 \nu \cdot x - \frac{1}{2} |\nabla_{\partial\Omega} u|^2 \nu \cdot x  - |x| |(\nu \cdot \nabla) u|  |\nabla_{\partial\Omega} u_{bd}|.
\]
Using the estimate $\nu \cdot x \geq -C|x|^2 \geq - Cr^2$ (see footnote \ref{Foot:nx}) in the first term, $|\nu \cdot x| \leq r$ and $|x| \leq r$ in the other terms, we get
\[
(\nu \cdot \nabla) u  \cdot (x \cdot \nabla) u - \frac{1}{2} |\nabla u|^2 \nu \cdot x
	\geq - Cr^2  |(\nu \cdot \nabla) u|^2 - \frac{1}{2} r |\nabla_{\partial\Omega} u|^2   - r |(\nu \cdot \nabla) u|  |\nabla_{\partial\Omega} u_{bd}|
\]
on $\partial \Omega \cap B_r$.
Altogether, we thus have:
\begin{align*}
\varphi'(r)
	&\geq -	\frac{1}{r^{M-2}} \int_{\partial \Omega \cap B_r}  \Big[ |(\nu \cdot \nabla) u|  |\nabla_{\partial\Omega} u_{bd}| + \frac{1}{2} |\nabla_{\partial\Omega} u_{bd}|^2  + Cr |(\nu \cdot \nabla) u|^2\Big]\,d\sigma\\
	&\geq - C r - \frac{C}{r^{M-3}}\int_{\partial \Omega \cap B_r} |(\nu \cdot \nabla) u|^2\,d\sigma - \frac{C}{r^{\frac{M-3}{2}}} \Big(\int_{\partial \Omega \cap B_r} |(\nu \cdot \nabla) u|^2\,d\sigma\Big)^{1/2},
\end{align*}

where $C>0$ depends on $\|\nabla_{\partial\Omega} u_{bd}\|_{L^\infty(\partial \Omega)}$ and $\partial\Omega$ (in particular on an upper bound of $\frac{|\partial\Omega \cap B_r|}{r^{M-1}}$ for $0 < r < r_0/2$). Appealing to the boundary Rellich estimate in Corollary \ref{Cor:BRE}, we have for every $0<r < r_0/2$:
$$\frac{1}{r^{M-3}}\int_{\partial \Omega \cap B_r} |(\nu \cdot \nabla) u|^2\,d\sigma\leq C\bigg(r^2+\varphi(2r)\bigg)$$ 
for a constant $C>0$ depending on $\partial \Omega$ and $\|\nabla_{\partial\Omega} u_{bd}\|_{L^\infty(\partial \Omega)}$. Combining the above estimates, the inequality \eqref{ineg22} follows. The last part of point 2. follows from the next lemma. 
\end{proof}

\begin{lemma}\label{lemma:supvarphi}
 Let $r_0 > 0$ and $\varphi: (0,r_0) \rightarrow [0,\infty)$ be a differentiable function satisfying, for some $C > 0$,
 \[
 \varphi'(r) \geq - Cr - C \varphi(2r) - C\sqrt{\varphi(2r)}, \quad \forall r\in (0, r_0/2).
 \]
 Then, for any $0 < r < r_0$,
 \begin{equation}
\sup_{(0,r]} \Big(1 + \sqrt{ \varphi(s) + \frac{1}{2}Cs^2}\Big) \leq e^{\frac{1}{4}Cr}\sup_{[r/2,r]} \Big(1 + \sqrt{ \varphi(s) + \frac{1}{2}Cs^2}\Big).
	\label{Eq:pDyaP}
 \end{equation}
 Consequently, for any $\delta \in (0,1)$, there exists $r_* > 0$ such that if $0 < r < r_*$ and $\varphi  \leq \delta$ in $[r/2,r]$, then $\varphi \leq 2\delta$ in $(0,r]$.
 
 \end{lemma}
 
\begin{proof}
 Note that once \eqref{Eq:pDyaP} is shown, the last assertion is obtained by simply selecting $r_* > 0$ sufficiently small so that
\[
e^{\frac{1}{4}Cr_*} \Big[1 + \big(\delta + \frac{1}{2} C r_*^2\big)^{1/2}\Big] - 1 \leq \sqrt{2\delta}.
\]
This is possible as the limit as $r_* \rightarrow 0$ of the left side is equal to $\sqrt{\delta} < \sqrt{2\delta}$.

Fix some $0 < r \leq r_0$. For $k\geq 0$, let $r_k = r 2^{-k}$,
\[
t_k = \sup_{(r_{k+1},r_k)} \varphi \quad \text{ and } \quad T_k = \sup_{(r_{k+1},r_k)} \Big[1 +  \big( \varphi(s) + \frac{1}{2}Cs^2\big)^{1/2}\Big].
\]
To prove \eqref{Eq:pDyaP}, we only need to show that
\begin{equation}
T_{k+1} \leq e^{\frac{1}{2}C r_{k+2}}T_k, \qquad k \geq 0.
	\label{Eq:pIter}
\end{equation}
(Once this is done, a simple induction gives
\[
T_{k+1} \leq e^{\frac{1}{2} C \sum_{k\geq 0} r_{k+2}} T_0 = e^{\frac{1}{4}Cr} T_0,  \hbox{ for all } k\geq 0, 
 \]
which yields the conclusion.) To prove \eqref{Eq:pIter}, consider for every $k \geq 0$ the open set:
\[
J_{k+1} = \{s \in (r_{k+2}, r_{k+1}): \varphi(2s) < \varphi(s)\} \subset (r_{k+2}, r_{k+1}).
\]
If $J_{k+1}$ is empty, we have $\varphi(s) \leq \varphi(2s)$ for $s \in [r_{k+2},r_{k+1}]$, which implies that $T_{k+1} \leq T_k$ and \eqref{Eq:pIter} follows. 

Suppose in the rest of the proof that $J_{k+1}$ is non-empty. For $s 
\in [r_{k+2},r_{k+1}] \setminus J_{k+1}$, we have $\varphi(s) \leq \varphi(2s)$ and so
\begin{equation}
\sup_{[r_{k+2},r_{k+1}] \setminus J_{k+1}} \Big[1 + \big( \varphi(s) + \frac{1}{2}Cs^2\big)^{1/2}\Big] \leq T_k.
	\label{Eq:pexJ}
\end{equation}
For the estimate in $J_{k+1}$, we use that $\varphi(2s) < \varphi(s)$ for $r\in J_{k+1}$ which combined with our hypothesis yields
$$\varphi'(s) \geq - C s - C \varphi(s) - C\sqrt{\varphi(s)}\quad \textrm{for every } \quad s\in J_{k+1}.$$
This implies 
\begin{equation}
\frac{d}{ds} \Big\{e^{\frac{1}{2}Cs}\Big[1 + \big(\varphi(s) + \frac{1}{2} Cs^2\big)^{1/2}\Big] \Big\} \geq 0 \quad \textrm{in } \quad J_{k+1}.
	\label{Eq:p1}
\end{equation}
Since $J_{k+1}$ is open, there exist at most countably many disjoint intervals $(a_{k,j},b_{k,j}) \subset (r_{k+2}, r_{k+1})$ such that $J_{k+1} = \cup_j (a_{k,j},b_{k,j})$. For each $j$, note that
\[
\varphi(b_{k,j}) = \begin{cases}
	\varphi(r_{k+1}) & \text{ if } b_{k,j} = r_{k+1},\\
	\varphi(2b_{k,j}) & \text{ if } b_{k,j} < r_{k+1},
	\end{cases}
\]
which implies
\begin{equation}
\varphi(b_{k,j}) \leq t_k.
	\label{Eq:p2}
\end{equation}
From \eqref{Eq:p1} and \eqref{Eq:p2}, we deduce that
\[
 \sup_{(a_{k,j},b_{k,j})}  \Big[1 + \big( \varphi(s) + \frac{1}{2}Cs^2\big)^{1/2}\Big] 
 	\leq e^{\frac{1}{2}C(b_{k,j} - a_{k,j})} \Big[1 + \big(t_k + \frac{1}{2}Cb_{k,j}^2\big)^{1/2}\Big] 
\leq e^{ \frac{1}{2}C r_{k+2}} T_k.
\]
Hence
\begin{equation}
 \sup_{J_{k+1}}  \Big[1 + \big( \varphi(s) + \frac{1}{2}Cs^2\big)^{1/2}\Big] 
 		\leq e^{ \frac{1}{2}C r_{k+2}} T_k.
	\label{Eq:pinJ}
\end{equation}
Estimate \eqref{Eq:pIter} is readily seen from \eqref{Eq:pexJ} and \eqref{Eq:pinJ}.
\end{proof}

%%%%%%%%%%%%%

\section{Uniform gradient estimate}\label{Sec:GradEst}

The main result of this section is the uniform gradient estimate for solutions $u_\eps$ to \eqref{eq:eps}
converging as in \eqref{ass:convH1} to a stationary harmonic map $u_0$ in the region where $u_0$ is continuous.

\begin{proposition} 
\label{pro:main}
Let $W \in C^2((-\infty,1])$ be such that \eqref{ass:Wstructure} and \eqref{ass:convex}  hold, $\Omega\subset \R^M$ be a $C^2$ bounded domain and 
$u_{bd}\in C^{1, \beta}(\partial \Omega, \Ss^{N-1})$ for some $\beta\in (0,1)$.
Let $(u_\eps)_{\eps\to 0}\subset H^1\cap L^\infty(\Omega, \R^N)$ be a family of solutions to \eqref{eq:eps} in $\Omega$ that has the common boundary value $u_{bd}$ and satisfies \eqref{ass:convH1} for a stationary harmonic map $u_0\in H^1(\Omega, \Ss^{N-1})$. Suppose $u_0$ is continuous in a neighborhood of a compact set $K\subset\bar\Omega$. Then there exists $C>0$ depending  only on $M$, $N$, $W$, $\Omega$, $K$, $u_{bd}$ and $u_0$ such that 
\[
\sup_{K}|\nabla u_{\eps}|\le C \quad \forall~ \eps > 0 . 
\]

\end{proposition}

The proof of Proposition \ref{pro:main} will be done in several steps:
\begin{enumerate}[(i)]
\item Continuity of $u_0$ implies small energy concentration: see Lemma \ref{lemma:sEc}.
\item Small energy concentration implies uniform convergence of $|u_\eps|$ to $1$: see Lemma \ref{lemma:unifconvmodule}.
\item Uniform convergence of $|u_\eps|$ implies Bochner's inequality: see Lemma \ref{lemma:Bochner}.
\item A pointwise bound of the gradient $\nabla u_\eps$ at the boundary in terms of energy concentration: see Lemma \ref{lem:bdryestim}.
\item Bochner's inequality and small energy concentration implies a uniform gradient bound: see Proposition \ref{prop:smallestim}.
\end{enumerate}
We note that steps (i) and (ii) do not require \eqref{ass:convex}, and that, in all steps except step (ii), one can replace \eqref{ass:Wstructure} by the slightly weaker condition \eqref{ass:Wsweak}.

The underlying ideas behind the above estimates have appeared many times in the literature. We refer the reader to e.g. \cite{GM79, HW75, Wiegner75} (for general nonlinear elliptic systems), \cite{su, SU-Bdry, schoen84} (for harmonic maps) and \cite{BBH93, BBO01, CLR18, LinRiv2, MZ10, NZ13} (for Ginzburg-Landau type problems). The novelty here concerns the very general assumptions on the potential $W$.

Let us begin by using the monotonicity formula in Proposition \ref{prop:monot} to show that the continuity of $u_0$ implies small energy concentration for $u_0$ (see Lemma \ref{lemma:u0sEc}) and subsequently for $u_\vareps$ (see Lemma \ref{lemma:sEc}).
\begin{lemma}\label{lemma:u0sEc}
Let $\Omega \subset \RR^M$ be a $C^2$ bounded domain and $u_0 \in H^1(\Omega, \mathbb{S}^{N-1})$ be a weakly harmonic map which is continuous on neighborhood of a compact set $K \subset \bar\Omega$. Suppose also that $u_0|_{\partial \Omega} = u_{bd}\in C^1(\partial\Omega, \mathbb{S}^{N-1})$. Then
\begin{multline}
\forall~ \eta > 0, \exists~ r_3 > 0   \text{ such that }
  \quad \frac{1}{r^{M-2}} \int_{\Omega \cap B_r(x_0)} |\nabla u_0|^2\,dx < \eta \quad \forall x_0 \in K, r \in (0,r_3].
  	\label{Eq:HMsEc}
\end{multline}
\end{lemma}

\begin{proof}
By hypothesis, $u_0$ is continuous in a larger compact set $K' \subset \bar\Omega$ of the form
\[
K' = (\bar B_{r_1}(0) + K) \cap \bar\Omega = \{x + z \in \bar \Omega: |x| \leq r_1, z \in K\}.
\]
In the proof, $C$ denotes a generic positive constant that depends only on $M$, $\Omega$, $u_{bd}$ and the modulus of continuity of $u_0$ on $K'$.  Since $u_0$ is uniformly continuous on $K'$, for every $\mu \in (0,1/8)$ there exists $r_2 \in (0,r_1 )$ such that
\begin{equation}
|u_0(x) - u_0(y)| \leq \mu \quad \forall~x, y \in K', |x - y| \leq  r_2.
	\label{Eq:u0UC}
\end{equation}

Fix some $\eta > 0$. In the following computation, we assume $x_0 \in K$, $r > 0$.

Case 1: $B_{2r}(x_0) \subset \Omega$ (that is $r \leq \mathrm{dist}(x_0,\partial\Omega)/2$). 

Pick $\chi \in C_c^\infty(B_{2r}(x_0))$ such that $0 \leq \chi \leq 1$, $|\nabla \chi| \leq \frac{2}{r}$ in $B_{2r}(x_0)$, and $\chi \equiv 1$ in $B_r(x_0)$. Testing the harmonic map equation against $(u_0 - u_0(x_0))\chi^2$, we obtain
\begin{align}
\int_{B_{2r}(x_0)} |\nabla u_0|^2 \chi^2 \,dx
	&\leq   \int_{B_{2r}(x_0)} |u_0 - u_0(x_0)|  \Big(2\chi |\nabla u_0 ||\nabla \chi| +  |\nabla u_0|^2 \chi^2 \Big)\,dx\nonumber\\
	&\leq   \int_{B_{2r}(x_0)}|u_0 - u_0(x_0)|  \Big(   |\nabla \chi|^2   +  2  |\nabla u_0|^2 \chi^2 \Big)\,dx.
	\label{Eq:HMSEC-1}
\end{align}
Imposing that $0 < r \leq r_3 \leq r_2/2$, we may use \eqref{Eq:u0UC} to upper bound $|u_0 - u_0(x_0)|$ on the right hand side of \eqref{Eq:HMSEC-1} and obtain
\[
\int_{B_{2r}(x_0)} |\nabla u_0|^2 \chi^2 \,dx \leq \frac{\mu}{1 - 2\mu}\int_{B_{2r}(x_0)}     |\nabla \chi|^2 \,dx \leq C \mu r^{M-2}.
\]
Recalling that $\chi \equiv 1$ in $B_{r}(x_0)$, we have attained \eqref{Eq:HMsEc} provided $\mu \leq \frac{\eta}{C}$, $0 < r \leq r_3 \leq r_2/2$ and $B_{2r}(x_0) \subset \Omega$.

Case 2: $\partial\Omega \cap B_{2r}(x_0) \neq \emptyset$ (that is $r > \mathrm{dist}(x_0,\partial\Omega)/2$). 

Pick $\tilde x_0 \in \partial\Omega \cap B_{2r}(x_0)$ so that $\Omega \cap B_{2r}(x_0) \subset \Omega \cap B_{4r}(\tilde x_0)$. Without loss of generality, we may assume $\tilde x_0 = 0$. As in the proof of Corollary \ref{Cor:BRE}, we can assume for some $r_0 = r_0(\Omega)$ that $\partial\Omega \cap B_{r_0} = \{(x',G(x')) \in \RR^M: x' \in U_0 \subset \RR^{M-1}\}$ and $\Omega \cap B_{r_0} = \{(x',x_M) \in B_{r_0}: x_M \geq G(x')\}$. Define 
\[
v_0((x',x_M)) = u_{bd}(x',G(x')) \text{ for } (x',x_M) \in \Omega \cap B_{r_0}.
\]

Pick $\chi \in C_c^\infty(B_{8r} )$ such that $0 \leq \chi \leq 1$, $|\nabla \chi| \leq \frac{2}{r}$ in $B_{8r} $, and $\chi \equiv 1$ in $B_{4r}$.  Imposing that $0 < r \leq r_3 \leq r_0/8$ and testing the harmonic map equation against $(u_0 - v_0)\chi^2$, we obtain
\begin{align}
\int_{\Omega \cap B_{8r}} |\nabla u_0|^2 \chi^2 \,dx
	&\leq    \int_{\Omega \cap B_{8r}} \Big[ |\nabla u_0||\nabla v_0| \chi^2 + |u_0 - v_0|\Big(2\chi  |\nabla u_0 ||\nabla \chi| \,dx +   |\nabla u_0|^2 \chi^2\Big)\Big] \,dx\nonumber\\
	&\leq   \Big(\int_{\Omega \cap B_{8r}}  |\nabla u_0|^2 \chi^2\,dx \Big)^{1/2} \Big( \int_{\Omega \cap B_{8r} }  |\nabla v_0|^2 \chi^2\,dx\Big)^{1/2} \nonumber\\
		&\qquad +  \int_{\Omega \cap B_{8r}}   |u_0 - v_0|\Big(  |\nabla \chi|^2  + 2|\nabla u_0|^2 \chi^2\Big) \,dx.
	\label{Eq:HMSEC-2}
\end{align}
Using the fact that $u_{bd} \in C^1(\partial\Omega)$, we have
\begin{equation}
\int_{\Omega \cap B_{8r}}  |\nabla v_0|^2 \chi^2\,dx \leq C(\|u_{bd}\|_{C^1(\partial\Omega)}) r^M.
	\label{Eq:HMSEC-3}
\end{equation}
Next, note that \eqref{Eq:u0UC} implies 
\[
|u_0(x) - u_0(0)| \leq \mu \quad \forall~x \in \bar \Omega \cap B_{r_2},
\]
which further implies
\[
|v_0(x) - u_0(0)| = |u_0(x', G(x')) - u_0(0)| \leq \mu \quad \forall~x \in \bar \Omega \cap B_{r_2} 
\]
and so
\begin{equation}
|u_0(x) - v_0(x)| \leq 2\mu \quad \forall~x \in \bar \Omega \cap B_{r_2} .
	\label{Eq:HMSEC-4}
\end{equation}
Imposing further that $  r_3 \leq r_2/8$ and using \eqref{Eq:HMSEC-3} and \eqref{Eq:HMSEC-4} in \eqref{Eq:HMSEC-2}, we have
\[
(1 - 2\mu) \int_{\Omega \cap B_{8r}} |\nabla u_0|^2 \chi^2 \,dx \leq Cr^{\frac{M}{2}} \Big(\int_{\Omega \cap B_{8r}}  |\nabla u_0|^2 \chi^2\,dx \Big)^{1/2}
	+ C \mu r^{M-2},
\]
which implies
\[
(1 - 4\mu) \int_{\Omega \cap B_{8r}} |\nabla u_0|^2 \chi^2 \,dx
	 \leq \frac{Cr^M}{\mu} 
	+ C \mu r^{M-2} \leq C\mu r^{M-2}(1 + \frac{r_3^2}{\mu^2}).
\]
Imposing further that $r_3 \leq \mu$, we get 
\[
\int_{\Omega \cap B_{8r}} |\nabla u_0|^2 \chi^2 \,dx
	 \leq  C \mu r^{M-2} .
\]
Recalling that $\chi \equiv 1$ in $B_{4r} \supset \Omega \cap B_r(x_0)$, we have attained  \eqref{Eq:HMsEc} provided $\mu \leq \frac{\eta}{C}$, $0 < r \leq r_3 \leq \min\{r_0/8,r_2/8,\mu\}$ and $\Omega \cap x_0 \in \partial \Omega$. This completes the  proof.
\end{proof}

\begin{remark}
Lemma \ref{lemma:u0sEc} remains valid if $u_{bd} \in C^{0,\beta}(\partial\Omega)$ for some $\beta \in (0,1]$. One only needs to replace the function $v_0$ in the proof by the harmonic extension $\tilde v_0$ of $u_{bd}$ in $\Omega$ and appeal to suitable estimates for the Dirichlet integral of $\tilde v_0$ on intersections of $\Omega$ with balls (\eqref{Eq:HMSEC-3} is replaced by a similar estimate where the term $r^M$ on the right hand side is replaced by $r^{M-2 + 2\beta}$), and slightly adjust the choice of $r_3$ towards the end. We omit the details.
\end{remark}

\begin{lemma}\label{lemma:sEc} Let $W \in C^2((-\infty,1])$ satisfy \eqref{ass:Wsweak}, and $\Omega\subset \R^M$ be a $C^2$ bounded domain. Let  $(u_\eps)_{\eps\to 0}\subset H^1\cap L^\infty(\Omega, \R^N)$ be a family of critical points to $\mcF_\vareps$ in $\Omega$ that  satisfies \eqref{ass:convH1} for a stationary harmonic map $u_0\in H^1(\Omega, \Ss^{N-1})$. Suppose $u_0$ is continuous in a neighborhood of a compact set $K\subset\bar\Omega$. If $K \cap \partial\Omega \neq \emptyset$, suppose in addition that $u_0|_{\partial\Omega} \in C^{1,\beta}(\partial\Omega)$, $\beta \in (0,1)$. Then 
\begin{multline}
\forall~\eta > 0, \exists~\eps_0 > 0, r_0 > 0 \text{ such that }\\
\frac{1}{r^{M-2}} \int_{\Omega \cap B_r(x_0)}e_\eps(u_\eps)\,dx \leq \eta \quad \forall~x_0 \in K, r \in (0,r_0], \eps \in (0,\eps_0].
	\label{Eq:sEc}
\end{multline}
Here $e_\eps(u_\eps)$ is the energy density in \eqref{def:edensity}, and the constants $\eps_0$, $r_0$ depend only on $M$, $N$, $\Omega$, $K$, $\eta$, $u_0$.
\end{lemma}

\begin{proof} By hypothesis, the limit $u_0$ is continuous in a larger 
compact set $K'\subset \bar\Omega$ of the form 
\[
K'=(\bar B_{r_1}(0)+K) \cap \bar\Omega=\{x+z \in \bar\Omega :\, |x|\leq r_1, z\in K\}
\]
 for some  $r_1>0$.

Recall that by \eqref{ass:Wsweak}, \eqref{ass:convH1} and Lemma \ref{Lem:u0ConvSHM}, $u_\eps \rightarrow u_0$ in $H^1(\Omega)$ and
\be\label{ass:vanishingpot}
\lim_{\eps\to 0}\int_\Omega \frac{W(1-|u_{\eps}(x)|^2)}{\eps^2}\,dx
	=0.
\ee

For clarity of exposion, we treat the cases $K \cap \partial\Omega \neq \emptyset$ and $K \cap \partial\Omega = \emptyset$ separately.

\medskip
\noindent
Step 1: Consider the situation that $K \cap \partial\Omega \neq \emptyset$.

Fix some $\gamma \in (0,1)$ for the moment. Let $r_*$ be the constant from point 2. of  Proposition~\ref{prop:monot} with $\delta = \gamma$. Let $d_* > 0$ be a small constant such that the distance function $d(x) := \textrm{dist}(x,\partial\Omega)$ is $C^2$ in $\{0 < d(x) < d_*\}$ and every point $x \in \{0 < d(x) < d_*\}$ has a unique closest point on $\partial\Omega$, denoted by $\pi(x)$. 

We claim that there exists $0 < r_0 \leq \min\{r_*/2,d_*, r_1/8\}$ and $\eps_0 > 0$ such that
\begin{equation}
\frac{1}{r_0^{M-2}}\int_{\Omega \cap B_{2r_0}(y)}e_\eps(u_\eps) \,dx \leq  \gamma \quad \forall~ y \in K'', 0 < \eps \leq \eps_0,
	\label{Eq:sEgamma}
\end{equation}
where
\[
K''=(\bar B_{r_1/8}(0)+K) \cap \bar\Omega=\{x+z \in \bar\Omega :\, |x|\leq r_1/8, z\in K\}.
\]

Indeed, by Lemma \ref{lemma:u0sEc}, there exists a small $0 < r_0 \leq \min\{r_*/2,d_*, r_1/8\}$ such that
$$ 
\int_{\Omega \cap  B_{2r_0}(y))} |\nabla u_0|^2\,dx  \le
\frac{r_0^{M-2}\gamma}{12} \quad \forall~ y \in K''.
$$

 Since $u_\eps\to u_0$ in
$H^1(\Omega)$, we can find $\vareps_0 > 0$ such that 
$\|\nabla u_\eps-\nabla u_0\|_{L^2(\Omega)}^2 \leq \frac{r_0^{M-2}\gamma}{12}$ for all $0 < \vareps \leq \vareps_0$
yielding by the triangle inequality:
$$
\int_{\Omega \cap B_{2r_0}(y)} |\nabla u_\eps|^2\,dx 
	\leq \big(\|\nabla u_\eps-\nabla u_0\|_{L^2(\Omega)}+ \|\nabla u_0\|_{L^2(\Omega \cap B_{2r_0}(y))}\big)^2
	\leq \frac{r_0^{M-2}\gamma}{3} 
		\quad \forall~ y \in K'', 0 < \eps \leq \eps_0.
$$
Recalling \eqref{ass:vanishingpot} and after possibly slightly shrinking $\eps_0$, we arrive at \eqref{Eq:sEgamma}.

We now fix $\eta > 0$, $x_0 \in K$, $r \in (0,r_0]$ and proceed to prove \eqref{Eq:sEc} by considering three different cases and indicating how $\gamma \leq \frac{1}{C(M,\Omega)}\eta$ is chosen (which then fix $\vareps_0$ and $r_0$ as above) in each case.

Case 1: $d(x_0) > r_0$. In this case, estimate \eqref{Eq:sEc} follows from \eqref{Eq:sEgamma} (with $y = x_0$) and the monotonicity formula in point 1. of Proposition~\ref{prop:monot} provided $\gamma \leq  \eta$. 

Case 2: $r_0 \geq d(x_0) \geq r$. Let $\tilde x_0 = \pi(x_0) \in K''$. Then
\[
B_{d(x_0)}(x_0) \subset \Omega \cap B_{2d(x_0)}(\tilde x_0).
\]
Using \eqref{Eq:sEgamma} with $y = \tilde x_0 \in \partial\Omega$ and noting that $2r_0 \leq r_*$, we have by the monotonicity formula in point 2. of Proposition~\ref{prop:monot} that
\[
\frac{1}{s^{M-2}} \int_{\Omega \cap B_s(\tilde x_0)} e_\eps(u_\eps) \,dx \leq 2  \gamma \text{ for all } s \in [0,2r_0].
\]
In particular, since $2d(x_0) \leq 2r_0$,
\[
\frac{1}{(2d(x_0))^{M-2}} \int_{\Omega \cap B_{2d(x_0)}(\tilde x_0)} e_\eps(u_\eps) \,dx \leq 2  \gamma.
\]
Hence
\[
\frac{1}{d(x_0)^{M-2}} \int_{B_{d(x_0)}(x_0)} e_\eps(u_\eps) \,dx 
	\leq \frac{1}{d(x_0)^{M-2}} \int_{\Omega \cap B_{2d(x_0)}(\tilde x_0)} e_\eps(u_\eps) \,dx 
	\leq 2^{M-1}  \gamma.
\]
Estimate \eqref{Eq:sEc} then follows for every $r \leq d(x_0)$ from the above inequality and the monotonicity formula in point 1. of Proposition~\ref{prop:monot} provided $\gamma \leq 2^{-(M-1)}\eta$.

Case 3: $d(x_0) < r \leq r_0$. Let $\tilde x_0 = \pi(x_0) \in K''$ as in Case 2. Note that 
\[
\Omega \cap B_{r}(x_0) \subset \Omega \cap B_{2r}(\tilde x_0) \text{ and } \frac{|\Omega \cap B_{r}(x_0)|}{|\Omega \cap B_{2r}(\tilde x_0) |} \geq \frac{1}{C(M,\Omega)}.
\]
As in Case 2, since $2r \leq 2r_0 \leq r_*$, we have by \eqref{Eq:sEgamma} (with $y = \tilde x_0$) and the monotonicity formula in point 2. of Proposition~\ref{prop:monot} that
\[
\frac{1}{(2r)^{M-2}} \int_{\Omega \cap B_{2r}(\tilde x_0)} e_\eps(u_\eps) \,dx \leq 2  \gamma.
\]
This implies
\[
\frac{1}{ r^{M-2}} \int_{\Omega \cap B_{ r}(  x_0)} e_\eps(u_\eps) \,dx \leq C(M,\Omega)  \gamma,
\]
which proves \eqref{Eq:sEc} provided $\gamma \leq \frac{\eta}{C(M,\Omega)}$. 

\medskip
\noindent
Step 2: Consider the situation that $K \cap \partial\Omega = \emptyset$. Arguing as in Step 1, we may choose the numbers $0 < r_0 < \min\{\mathrm{dist}(K,\partial\Omega)/2, r_1/8\}$ such that \eqref{Eq:sEgamma} holds with $\gamma = \eta$. Estimate \eqref{Eq:sEc} then follows from \eqref{Eq:sEgamma} as in Case 1 of Step 1. The proof is completed.
\end{proof}

We next show that small energy concentration implies that $|u_\eps| \rightarrow 1$ uniformly.

\begin{lemma}\label{lemma:unifconvmodule} 
Let $W \in C^2((-\infty,1])$ satisfy \eqref{ass:Wstructure}, $\Omega\subset \R^M$ be a $C^2$ bounded domain, and
$u_{bd}\in C^{1, \beta}(\partial \Omega, \Ss^{N-1})$, $\beta \in (0,1)$. Let  $(u_\eps)_{\eps\to 0}\subset H^1\cap L^\infty(\Omega, \R^N)$ be a family of solutions to \eqref{eq:eps}  such that \eqref{ass:convH1} holds for a stationary harmonic map $u_0\in H^1(\Omega, \Ss^{N-1})$. Suppose that \eqref{Eq:sEc} holds for a compact set $K \subset \bar\Omega$. Then
$$
|u_\eps|\to 1 \quad \textrm{ uniformly in } K \textrm{ as } \eps\to 0.
$$
\end{lemma}
\begin{proof} By Theorem~\ref{thm:MPI}, we have
\[
|u_\eps| \leq 1 \text{ in } \Omega.
\]
In the proof below,  $C$ denotes a generic positive constant which depends only on the dimension $M$, $N$, $\Omega$, $K$, $u_{bd}$, $u_0$ and the behaviour of the function $W$ in $[0, 1]$. 

By \eqref{ass:Wstructure}, it suffices to show
\begin{equation}
0 \leq \alpha_\eps :=W(1-|u_\eps|^2) \to 0 \text{ uniformly in $K$ as } \eps \to 0.
	\label{Eq:alepu0bdry}
\end{equation}
We now argue that, in order to show \eqref{Eq:alepu0bdry}, we only need to show that there exist $C_1, C_2 > 0$ and $\bar \eps > 0$ such that, with $r_\eps(x_0) := \frac{\alpha_\eps(x_0)\vareps}{2C_1}$,
\begin{equation}
 \alpha_\eps^3(x_0)
 	\le 
\frac{C_2}{r_\eps(x_0)^{M-2}} \int_{\Omega \cap B_{r_\eps(x_0)}(x_0)}\frac{ \alpha_\vareps(x)}{\eps^2}\,dx
	\quad \forall~x_0\in K, 0 < \eps \leq \bar \eps.
\label{sabXbdry}
 \end{equation}
Indeed, for any given $\eta > 0$, and with $\eps_0$ and $r_0$ as in \eqref{Eq:sEc}, we have for $0 < \eps \leq \min\big\{\eps_0, \bar \eps, \frac{2C_1r_0}{\max_{[0,1]}W}\big\}$ and $x_0 \in K$ that $r_\eps(x_0) \leq r_0$ and
\[
\alpha_\eps^3(x_0)
 	\stackrel{\eqref{sabXbdry}}{\le} \frac{C_2}{r_\eps(x_0)^{M-2}} \int_{\Omega \cap B_{r_\eps(x_0)}(x_0)}\frac{ \alpha_\vareps(x)}{\eps^2}\,dx 
 	\stackrel{\eqref{Eq:sEc}}{\le} C_2\eta.
\]
Since $\eta$ is arbitrary, this implies \eqref{Eq:alepu0bdry}.

Before proving \eqref{sabXbdry}, we claim that  
\begin{equation}
\label{eq:bulk3} \|\nabla u_\eps\|_{L^\infty(\Omega)}^2 \leq
\frac{C}{\eps^2}.
\end{equation}
Indeed, we decompose $u_\eps=v_\eps+w$ where $w$ is the harmonic function that is equal to $u_{bd}$ on $\partial\Omega$ and $v_\eps \in H_0^1(\Omega)$ satisfies $-\Delta v_\eps = \frac{1}{\eps^2} W'(1 - |u_\eps|^2)u_\eps$ in $\Omega$. Then, by elliptic estimates,
\[
\|w\|_{L^\infty(\Omega)} + \|\nabla w\|_{L^\infty(\Omega)}\le C,
\]
and, by \cite[Lemma $A.2$]{BBH93},
\[
\|\nabla v_\eps\|_{L^\infty(\Omega)}^2 \le C\|\Delta v_\eps\|_{L^\infty(\Omega)}\|v_\eps\|_{L^\infty(\Omega)} \leq \frac{C}{\eps^2}.
\]
The claim \eqref{eq:bulk3} follows.

We now prove \eqref{sabXbdry}. As $W\in C^1$ and $|u_\eps|\leq 1$ in $\Omega$, we have by \eqref{eq:bulk3}: 
\[
|\alpha_\eps (x_0) -  \alpha_\eps(x) |  
	\le   C |u_\eps(x)-u_\eps(x_0)|
	\le 	\frac{C_1}{\eps}|x-x_0|
	 \quad \forall~x\in \bar\Omega.
\]
Recalling $r_\vareps(x_0) = \frac{\alpha_\eps(x_0) \eps}{2C_1}$, we then have $x_0 \in K$ that 
\[
 \alpha_\eps(x_0)\leq 2 \alpha_\eps(x)  \quad \forall
x\in \overline{\Omega \cap B_{r_\eps(x_0)}(x_0)}.
\]
Dividing the above inequality by $2r_\eps(x_0)^{M-2}\eps^2$,
 integrating over $\Omega \cap B_{r_\eps(x_0)}(x_0)$, and noting that
\[
|\Omega \cap B_{r_\eps(x_0)}(x_0)| \geq \frac{1}{C}r_\eps(x_0)^M,
\]
we obtain  
\[
 \alpha_\eps^3(x_0) = \frac{4C_1^2r_\eps(x_0)^2}{\eps^2}\alpha_\eps(x_0)\le 
\frac{C_2}{r_\eps(x_0)^{M-2}} \int_{\Omega \cap B_{r_\eps(x_0)}(x_0)}\frac{ \alpha_\vareps(x)}{\eps^2}\,dx,
\]
which gives \eqref{sabXbdry} and concludes the proof.
\end{proof}

We note the following interior version of Lemma \ref{lemma:unifconvmodule}, which is not needed in the present paper.

\begin{lemma}\label{lemma:unifconvmoduleint} Let $W \in C^2((-\infty,1])$ satisfy \eqref{ass:Wstructure}, $\Omega\subset \R^M$ be a domain, and $(u_\eps)_{\eps\to 0}\subset H^1\cap L^\infty(\Omega, \R^N)$ be a family of critical points of $\mcF_\vareps$ in $\Omega$ that is uniformly bounded in $L^\infty(\Omega)$ and satisfies \eqref{ass:convH1} for a stationary harmonic map $u_0\in H^1(\Omega, \Ss^{N-1})$. Suppose that \eqref{Eq:sEc} holds on a neighborhood of a compact set $K \subset \Omega$. Then
$$
|u_\eps|\to 1 \quad  \text{ uniformly on }K \text{ as } \eps\to 0.
$$
\end{lemma}

\begin{proof} The proof is similar to that of Lemma \ref{lemma:unifconvmodule}, except that instead of \eqref{eq:bulk3}, one uses the estimate
 \be
 \label{estim33} 
 \|\nabla u_\eps\|_{L^\infty(K')}^2\le C \left(\|\Delta u_\eps\|_{L^\infty(\Omega)}+\frac{1}{\textrm{dist}^2(K',\partial\Omega)}\right)
\leq \frac{C}{\eps^2} \quad \forall~0 < \eps   \leq  \mathrm{dist}(K',\partial\Omega), 
\ee
which follows from \cite[Lemma A.1]{BBH93}. We omit the details.
\end{proof}

Thanks to Lemma \ref{lemma:unifconvmodule}, we can assume that $\big|1-|u_\eps|\big|$ is uniformly small in the region where the limit $u_0$ is continuous. In that region, the following Bochner's inequality holds for the density energy $e_\eps(u_\eps)$.

\begin{lemma}[Bochner's inequality]\label{lemma:Bochner}
Let $W\in C^2((-\infty, 1])$ be such that \eqref{ass:Wsweak} and \eqref{ass:convex} hold with some $\kappa >0$. Let $0 < \delta \leq \min\{1/4, \kappa/3\}$. If $u_\eps\in H^1\cap L^\infty(\Omega)$ is a critical point of $\mcF_\vareps$ with $\eps > 0$ in a domain $\Omega\subset \R^M$ such that $\big|1-|u_\eps|\big|\leq\delta$ in a subdomain $\omega \subset \Omega$, then
\be\label{est:Bochner}
-\Delta [e_\eps(u_\eps)]+|\nabla^2 u_\eps|^2 \le  2|\nabla u_\eps|^4\le 8e_\eps(u_\eps)^2 \quad \textrm{ in } \, \omega,
\ee
where $e_\eps(u_\eps)$ is the density energy in \eqref{def:edensity}. 
\end{lemma}
 
\begin{proof} For simplicity, we drop the index $\eps$ in the following. As $u\in H^1\cap L^\infty(\Omega)$ is a solution to \eqref{eq:eps}, elliptic theories give that $u\in W^{3,p}_{\loc}\cap C^2(\Omega)$ for every $p<\infty$. We claim that:
$$
\Delta [e(u)]=|\nabla^2 u|^2+|\Delta u|^2-\frac{2|\nabla u|^2}{\eps^2}W'(1-|u|^2)+\frac{|\nabla |u|^2|^2W''(1-|u|^2)}{\eps^2}
\quad \textrm{a.e.  in } \, \Omega;
$$
in particular, $\Delta [e(u)]\in C^0(\Omega)$ as the above RHS is continuous in $\Omega$.  
Indeed, we compute firstly a.e. in $\Omega$:
$$
\frac12 \Delta [|\nabla u|^2]=\frac12 \sum_{i,j,k} \nabla_k^2 [(\nabla_i u_j)^2]= \sum_{i,j,k} \big[(\nabla_k \nabla_i u_j)^2+(\nabla_k^2 \nabla_i u_j) \nabla_i u_j\big] =|\nabla^2 u|^2+ \sum_{i,j} ( \nabla_i \Delta u_j) \nabla_i u_j.
$$
To compute the last term above, we use for every $i,j$:
$$
-\nabla_i \Delta u_j=\frac{\nabla_i u_j}{\eps^2} W'(1-|u|^2) - \sum_k\frac{2 u_j u_k \nabla_i u_k }{\eps^2} W''(1-|u|^2)
\quad \textrm{a.e. in} \, \Omega,$$
yielding 
$$
\sum_{i,j}(\nabla_i \Delta u_j) \nabla_i u_j = -\frac{|\nabla u|^2}{\eps^2} W'(1-|u|^2) + \frac{|\nabla |u|^2|^2 }{2\eps^2} W''(1-|u|^2) \quad \textrm{a.e. in} \, \Omega.
$$
Secondly, as $W(1-|u|^2)\in C^2(\Omega)$, we compute pointwise in $\Omega$: 
\begin{align*}
\frac{1}{2\eps^2}\sum_k \nabla_k^2 \big[ W(1-|u|^2)\big] &=\frac{1}{2\eps^2} \sum_k \nabla_k [-2 u\cdot \nabla_k u W'(1-|u|^2)] \\
&=-\frac{1}{\eps^2} (u\cdot \Delta u +|\nabla u|^2) W'(1-|u|^2) + \frac{|\nabla |u|^2|^2}{ 2\eps^2} W''(1-|u|^2).
\end{align*}
Summing up, we obtain a.e. in $\Omega$:
\begin{align*}
\Delta [e(u)] &=  |\nabla^2 u|^2 -\frac{2 |\nabla u|^2}{\eps^2} W'(1-|u|^2) -\frac{u\cdot \Delta u }{\eps^2} W'(1-|u|^2) + \frac{|\nabla |u|^2|^2 }{\eps^2} W''(1-|u|^2) \\
&= |\nabla^2 u|^2 -\frac{2 |\nabla u|^2}{\eps^2} W'(1-|u|^2)+|\Delta u|^2  + \frac{|\nabla |u|^2|^2 }{\eps^2} W''(1-|u|^2)
\end{align*}
which proves the claim.

We now turn our attention to $\omega$ where $\big|1-|u|\big|\le \delta \leq \frac{1}{4}$. On this subdomain, we have  $|u|^2 \geq (1-\delta)^2\geq \frac12$ and so, by \eqref{eq:eps}, 
$$
\frac{1}{2\eps^4}\left(W'(1-|u|^2)\right)^2\leq |\Delta u|^2 \quad \text{ in } \omega.
$$
Using that
$$
\frac{2|\nabla u|^2}{\eps^2}W'(1-|u|^2)\le\frac{1}{2\eps^4}(W'(1-|u|^2))^2+ 2|\nabla u|^4 ,
$$ 
we conclude
\begin{align*}
-\Delta [e(u)] +  |\nabla^2 u|^2 &=\frac{2 |\nabla u|^2}{\eps^2} W'(1-|u|^2)  -|\Delta u|^2 - \frac{|\nabla |u|^2|^2 }{\eps^2} W''(1-|u|^2) \\
&\leq  2|\nabla u|^4 - \frac{|\nabla |u|^2|^2 }{\eps^2} W''(1-|u|^2) \leq 2|\nabla u|^4 \quad \text{ in } \omega.
\end{align*}
Above, we used \eqref{ass:convex} thanks to $|1-|u|^2| \leq (2+\delta)\delta< 3 \delta \leq \kappa$. Estimate \eqref{est:Bochner} follows from the nonnegativity of $W$ in \eqref{ass:Wsweak}.
\end{proof}

Furthermore, in the regime where $1-|u_\eps|$ is non-negative (thanks to Theorem \ref{thm:MPI}) and small enough (thanks to Lemma \ref{lemma:unifconvmodule}), the Laplacian of the modulus $\rho_\eps=|u_\eps|$ and of the director $\ds n_\eps=\frac{u_\eps}{\rho_\eps}$ of $u_\eps$ are controlled by the energy density $e_\eps(u_\eps)$. These estimates do not use Bochner's inequality.

\begin{lemma}\label{lemma:Laplaciantannorm}
Let $W\in C^2((-\infty, 1])$ be such that \eqref{ass:Wsweak} and \eqref{ass:convex} hold  with some $\kappa >0$. Let $0 < \delta \leq \min\{1/4, \kappa/3\}$. Let $u_\eps\in H^1\cap L^\infty(\Omega)$ be a critical point of $\mcF_\vareps$ with $\eps > 0$ in a domain $\Omega\subset \R^M$ such that $0\leq 1-|u_\eps|\leq\delta$ in a subdomain $\omega \subset \Omega$. Then, within the decomposition $u_\eps:=\rho_\eps n_\eps$ where $\rho_\eps:=|u_\eps|$ and $n_\eps\in\mathbb{S}^{N-1}$,  we have
$$
 |\Delta n_\eps|\le 8e_\eps(u_\eps)
\quad \textrm{ and }
\quad
 \Delta \rho_\eps\le 4e_\eps(u_\eps) \quad \textrm{ in } \, \omega.
$$

\end{lemma}

\begin{proof} For simplicity, we drop the index $\eps$ in the following. We know that $u\in W^{3,p}_{\loc}\cap C^2(\Omega)$ for every $p<\infty$. As $|u|\geq \frac 12$ in $\omega$, we deduce that $\rho = |u|$ and $n=\frac{u}{\rho}$ have the same regularity in $\omega$. We compute in $\omega$: $n \cdot \Delta n +|\nabla n|^2 =0$ (as $\Delta (|n|^2)=0$), $n\cdot \nabla_k n=0$ for every $1\leq k\leq M$ and
$$
\Delta u=\Delta (\rho n)=n\Delta \rho + \rho \Delta n + 2 (\nabla \rho \cdot \nabla) n.
$$
These identities combines with \eqref{eq:eps} multiplied by $n$ yield
\be\label{eq:rho}
\Delta \rho+\frac{\rho}{\eps^2}W'(1-\rho^2)=\rho|\nabla n|^2 \quad \textrm{in} \, \omega.
\ee
Also, multiplying \eqref{eq:eps} by any orthogonal vector $m \perp n(x)$ at a point $x\in \omega$, we obtain
$\ds 
0=m\cdot \Delta u(x)=m \big(\rho(x) \Delta n(x) +  2 \nabla \rho(x)\cdot \nabla n(x)\big)$.
This yields
\begin{equation}
\rho \Delta n +  2 (\nabla \rho \cdot \nabla) n =\lambda n= -\rho |\nabla n|^2 n \quad \textrm{in} \, \omega.
	\label{eq:n}
\end{equation}

As $W$ is non-negative by \eqref{ass:Wsweak}, we have
\begin{equation}
2e(u) \geq |\nabla u|^2 = |\nabla \rho|^2 + \rho^2|\nabla n|^2.
	\label{eq:eu>du}
\end{equation}
From \eqref{eq:n}, \eqref{eq:eu>du} and the fact that $\rho \geq \frac{3}{4}$ in $\omega$, it is readily seen that
\[
|\Delta n| \leq  \frac{1}{\rho^2} |\nabla \rho|^2 + 2|\nabla n|^2   \leq \frac{4}{\rho^2} e(u) \leq  8 e(u) \quad \textrm{in} \, \omega.
\]
On the other hand, by the assumptions on $\rho$, we have $0\leq 1-\rho^2\leq 3\delta \leq \kappa$, and thus, by \eqref{est:Wincreas} (which is a consequence of \eqref{ass:Wsweak} and \eqref{ass:convex}), $W'(1-\rho^2)\geq 0$ in $\omega$. Using this as well as \eqref{eq:eu>du} and the fact that $\rho \geq \frac{3}{4}$ in $\omega$ in \eqref{eq:rho}, we get
$$
\Delta \rho\le \rho|\nabla n|^2 \leq \frac{2}{\rho} e(u)   \leq 4e(u) \quad \textrm{ in } \, \omega.
$$ 
\end{proof}

As a consequence of Lemma \ref{lemma:Laplaciantannorm}, we obtain the following boundary gradient estimates. (See \cite[Lemma 3.3]{CLR18} for a closely related result where the boundary map is assumed to be twice differentiable.):

\begin{lemma}\label{lem:bdryestim} Let $W \in C^2((-\infty,1])$ be such that \eqref{ass:Wsweak} and \eqref{ass:convex} hold with some $\kappa >0$. Let $0 < \delta \leq \min\{1/4, \kappa/3\}$. Let $\Omega\subset \R^M$ be a $C^2$ bounded domain and 
$u_{bd}\in C^{1, \beta}(\partial \Omega, \Ss^{N-1})$ for some $\beta\in (0,1)$. 
Then for every $p>M$, 
there exist $C, r_0>0$ (depending on $M$, $\beta$, $p$  and $\Omega$) such that for any $x_0\in\bar\Omega$, $r\in (0, r_0)$,  $\eps > 0$, and any solution $u_\eps\in H^1\cap L^\infty(\Omega, \R^N)$ to \eqref{eq:eps}, if
\begin{equation*}
1-|u_\eps|\leq\delta \mbox{ in } B_r(x_0)\cap\Omega,
\end{equation*}
then
it holds 
\begin{eqnarray*}
\sup_{B_{r/2}(x_0)\cap \partial\Omega} |\nabla u_\eps| &\leq& C \left( r^{{1-\frac Mp}}\norm{e_{ \eps}(u_\eps)}_{L^p(B_r(x_0)\cap\Omega)}  \right.
+ r^{-\frac{M}2}\norm{\nabla u_\eps}_{L^2(B_r(x_0)\cap\Omega)} \\
&\quad& +\left. \|\nabla_\tau u_{bd}\|_{L^\infty(B_r(x_0)\cap\partial\Omega)}+   r^\beta [\nabla_{\tau} u_{bd}]_{C^\beta(B_r(x_0)\cap\partial \Omega)}+\frac\delta r \right),
\end{eqnarray*}
where $\ds [\nabla_{\tau} u_{bd}]_{C^\beta(\partial \Omega)}=\sup_{x\neq y, x,y\in \partial \Omega}\frac{|\nabla_{\tau} u_{bd}(x)-\nabla_{\tau} u_{bd}(y)|}{|x-y|^\beta}$ and $\nabla_{\tau} u_{bd}$ denotes the gradient vector of $u_{bd}$ as a function defined on the submanifold $\partial\Omega$ of $\RR^M$.
\end{lemma}

\begin{remark}
The constant $C$ does not depend on the boundary data $u_{bd}$ nor the function $W$.
\end{remark}

\begin{proof} As before, we drop the index $\eps$ and we write $B_r$ instead of $B_r(x_0)$ in the following.
Recall that $u\in C^2(\Omega)\cap  C^{1,\beta}(\Omega)$ (see Remark \ref{rem:regu}) and $|u|\leq 1$ in $\Omega$ by Theorem \ref{thm:MPI}. 

\bigskip

\nd {\it Case 1. $x_0\in \partial \Omega$}. We abbreviate $B_r(x_0)$ to $B_r$. 
We use the decomposition $u=\rho n$ with $\rho=|u|$ and $\ds n=\frac{u}{\rho}$. As $|\nabla n|\leq 2 |\nabla u|$ in $B_r\cap \Omega$ and $n=u_{bd}$ on $B_r\cap \partial \Omega$, we use Lemma \ref{lemma:Laplaciantannorm} and Lemma \ref{lem:bdry-est} (see Appendix) 
to get 
\begin{multline}\label{estimgraduN}
\sup_{B_{r/2}\cap\Omega}\abs{\nabla n} \leq C\Big(r^{1- \frac{M}{p}} \norm{e(u)}_{L^p(B_r\cap\Omega)} + r^{-\frac{M}{2}} \norm{\nabla u}_{L^2(B_r\cap\Omega)}\\
	+  \norm{\nabla_\tau u_{bd}}_{L^\infty(B_r\cap \partial \Omega)}+ r^\beta [\nabla_\tau u_{bd}]_{C^\beta(B_r\cap \partial \Omega)}  \Big),
\end{multline}
where here and below $C$ denotes a constant depending only on $M$, $\beta$, $p$ and $\Omega$. It remains to estimate $\nabla \rho$ in $L^\infty(B_{r/2}\cap\partial \Omega)$. As $\rho=|u_{bd}|=1$ on $\partial \Omega$, in particular, $\nabla _\tau \rho=0$ on $\partial \Omega$, we only need to estimate $\frac{\partial \rho}{\partial \nu}$ in $L^\infty(B_{r/2}\cap \partial \Omega)$. For that, we denote $\tilde\rho=1-\rho \geq 0$ in $B_r\cap \Omega$; by Lemma \ref{lemma:Laplaciantannorm}, we have $-\Delta \tilde\rho\le 4e(u)$ in $B_r\cap \Omega$. Consider the solution\footnote{Such a solution exists and is unique because $e(u)\in L^p(B_r\cap \Omega)$, $p\geq 2$ and $\tilde\rho$ is Lipschitz on $\partial (B_r\cap \Omega)$.} $w\in H^1(B_1\cap\Omega)$ to the Dirichlet problem
\begin{equation*}
-\Delta w = 4e(u) \text{ in }B_r\cap\Omega,\qquad
w=\tilde\rho  \text{ on }\partial(B_r\cap\Omega).
\end{equation*}
By the maximum principle, $\tilde\rho\leq w$ in $B_r\cap\Omega$. As $w=\tilde\rho=0$ on $B_r\cap\partial\Omega$, Lemma \ref{lem:bdry-est} (see Appendix)  implies
\begin{equation*}
\sup_{B_{r/2}\cap\Omega}\abs{\nabla w}\leq C \left( r^{1 - \frac{M}{p}}\norm{e(u)}_{L^p(B_{3r/4}\cap\Omega)} + r^{-\frac{M}{2}} \norm{\nabla w}_{L^2(B_{3r/4}\cap\Omega)}\right).
\end{equation*}
To estimate the last term of the above inequality, one decomposes $w=w_1+w_2$, where $\Delta w_2=0$ and $-\Delta w_1=4 e(u)$  in $B_r\cap\Omega$, with $w_1$ vanishing on the full boundary $\partial(B_r\cap\Omega)$ and $w_2=\tilde \rho\in [0, \delta]$ on $\partial(B_r\cap\Omega)$. Then 
\begin{align*}
\norm{\nabla w_1}_{L^2(B_{3r/4}\cap\Omega)}& \leq  C r^{1 - \frac{M}{p} + \frac{M}{2}}\norm{e (u)}_{L^p(B_r\cap\Omega)},\\
\norm{\nabla w_2}_{L^2(B_{3r/4}\cap\Omega)}&\leq C r^{-1 + \frac{M}{2}} \norm{w_2}_{L^\infty(B_r\cap\Omega)}\leq C r^{-1 + \frac{M}{2}}\norm{w_2}_{L^\infty(\partial(B_r\cap\Omega))}\leq Cr^{-1 + \frac{M}{2}}\delta.
\end{align*}
Summing up, we deduce
\begin{equation*}
\sup_{B_{r/2}\cap\Omega}\abs{\nabla w}\leq C \left( r^{1 - \frac{M}{p}} \norm{e (u)}_{L^p(B_{r}\cap\Omega)} + \frac{\delta}{r} \right).
\end{equation*}
Since  $w=\tilde \rho=0$ on $B_r\cap\partial \Omega$, this entails
\begin{equation*}
0\leq \tilde\rho\leq w \leq C \left( r^{1 - \frac{M}{p}} \norm{e (u)}_{L^p(B_r\cap\Omega)}+\frac{\delta}{r}\right) \dist(\cdot,\partial\Omega) \quad \textrm{in } B_{r/2}\cap\partial \Omega,
\end{equation*}
which imply
\begin{equation}\label{estimgraduperp}
\sup_{B_{r/2}\cap\partial\Omega}\abs{\frac{\partial \rho}{\partial\nu}}
	=\sup_{B_{r/2}\cap\partial\Omega}\abs{\frac{\partial \tilde \rho}{\partial\nu}}
	\leq C \left( r^{1 - \frac{M}{p}}\norm{e(u)}_{L^p(B_r\cap\Omega)}+ \frac{\delta}{r}\right).
\end{equation}
The conclusion follows from \eqref{estimgraduN}-\eqref{estimgraduperp}.

\medskip

\nd {\it Case 2. $x_0\in \Omega$}. We first fix $r_0$ as follows. Since $\partial \Omega$ is $C^2$, by a similar argument as in the beginning of the proof of Corollary \ref{Cor:BRE}, there exists $r_0>0$ depending only on $\Omega$ such that for every $s \in (0, r_0]$ and any $x\in \Omega$ with $s> 2\dist(x,\partial \Omega)$, we have $\mathrm{diam}(B_{2s/3}(x)\cap \partial \Omega)\geq s/6$. 

If $r \leq 2\dist(x_0,\partial \Omega)$ (where $B_{r/2}(x_0)\cap \partial\Omega=\emptyset$), the conclusion is obvious. Thus, we may assume that $r > 2\dist(x_0,\partial \Omega)$. By the property of $r_0$, $B_{2r/3}(x_0)\cap \partial \Omega$ can be covered by a finite number (depending only on $\partial \Omega$) of  balls $\{B_{r/12}(y_j)\}_{j\in J}$ with $y_j\in \partial \Omega$ such that
$B_{r/6}(y_j)\subset B_{r}(x_0)$ for all $j\in J$. Applying Case 1 for the balls $B_{r/6}(y_j)$ for $j\in J$, up to eventually increasing the constant $C$, the desired conclusion follows in the ball $B_{r}(x_0)$. 
\end{proof}

Now we prove the $\eta$-regularity result for solutions $u_\eps$ to \eqref{eq:eps}.

\begin{prop}\label{prop:smallestim} 
Let $W \in C^2((-\infty,1])$ be such that \eqref{ass:Wsweak} and \eqref{ass:convex} hold  with some $\kappa >0$. Let $\Omega\subset \R^M$ be a $C^2$ bounded domain and 
$u_{bd}\in C^{1, \beta}(\partial \Omega, \Ss^{N-1})$ for some $\beta\in (0,1)$. 
Then there exist $\eta_0>0$, $\delta > 0$, $r_0 >0$ and $C > 0$ (depending only on $M$, $\beta$, $\kappa$ and $\Omega$) such that 
for every $\eps >0$, {$r\in (0,r_0)$}, $x_0\in\overline\Omega$ and any solution $u_\eps\in H^1\cap L^\infty(\Omega, \R^N)$ to \eqref{eq:eps} satisfying
\begin{equation*}
 1-|u_\eps| \leq\delta \mbox{ in } B_r(x_0)\cap\Omega,
\end{equation*}
and
\begin{equation*}
E:=\sup_{B_s(x)\subset B_{r}(x_0)} \; \frac1{s^{M-2}}\int_{\Omega\cap B_s(x)} e_\eps(u_\vareps) \leq \eta_0,
\end{equation*}
there holds
\begin{equation*}
\sup_{B_{r/2}(x_0)\cap\Omega}\, e_\eps(u_\eps) \leq C  \left( \frac{E + \delta^2}{r^2} + \|\nabla_\tau u_{bd}\|_{C^{0,\beta}(\partial\Omega)}^2 \right),
\end{equation*}
where $\|\cdot\|_{C^{0,\beta}(\partial\Omega)}=\|\cdot\|_{L^\infty(\partial\Omega)}+[\cdot]_{C^{\beta}(\partial\Omega)}$.
\end{prop}

\begin{remark}
The constants $\eta_0, \delta, r_0, C$ do not depend on the boundary data $u_{bd}$. They depend on the function $W$ only through the constant $\kappa$.
\end{remark}

\begin{proof} The proof is similar to the proof of $\eta$-regularity result for harmonic maps, see e.g.  \cite{schoen84}. 

Recall that $u\in C^2(\Omega)\cap  C^{1,\beta}(\Omega)$ (see Remark \ref{rem:regu}) and $|u|\leq 1$ in $\Omega$ by Theorem \ref{thm:MPI}. For $r>0$, consider the quantity
\begin{equation}\label{smallestim1}
\mathcal{M}^2 :=\sup_{0<s<r} \bigg((r-s)^2 \sup_{B_s(x_0) \cap\Omega}  e_\e (u_\eps)\bigg).
\end{equation}
We aim to prove for small $r>0$:
\begin{equation}
\mathcal{M}^2 \leq C  \left( E + \delta^2 + \|\nabla_\tau u_{bd}\|_{C^{0,\beta}(\partial\Omega)}^2 r^2  \right),
	\label{Eq:MGoal}
\end{equation}
for $C >0$ depending only on $M$, $\beta$, $W$ and $\Omega$. This implies the conclusion by taking $s=r/2$ in 
\eqref{smallestim1}.

Throughout the proof $\delta$ and $r_0$ are sufficiently small such that Lemma \ref{lem:bdryestim} as well as the Bochner-type inequality \eqref{est:Bochner} hold for $ u$ in $B_{r}(x_0) \cap \Omega$ for every $r\in (0, r_0)$. If the supremum over $s$ in \eqref{smallestim1} is attained by $s = r$, then $\mathcal{M} = 0$ and \eqref{Eq:MGoal} clearly holds. We may thus assume without loss of generality that $\mathcal{M}>0$ and then, as $u_\eps\in C^1(\bar \Omega)$, there exist $s_0\in [0,r)$  and $x_1\in \overline {B_{s_0}(x_0) \cap \Omega}$ such that
\begin{equation*}
\mathcal{M}^2 =(r-s_0)^2 \sup_{B_{s_0}(x_0) \cap\Omega}  e_\e (u_\eps)   = (r-s_0)^2  e_\e (u_\eps)(x_1)  .
\end{equation*}
With $s_1:=(r-s_0)/2 > 0$, it holds $\mathcal{M}^2 =4s_1^2   e_\e (u_\eps)(x_1) $ and 
\begin{align*}
\sup_{B_{s_1}(x_1)\cap\Omega}e_\e (u_\eps) 
	& \leq \sup_{B_{s_1+s_0}(x_0) \cap\Omega}e_\e (u_\eps)
	 \leq \frac{\mathcal{M}^2}{(r-s_1-s_0)^2}  	
	 =\frac{\mathcal{M}^2}{s_1^2} 
	 =4 e_\e (u_\eps)(x_1).
\end{align*}
In the following, we use the rescaling: for $V:=e_\e (u_\eps)(x_1)>0$ (as $\mathcal M>0$), we set $\widetilde\Omega:=\sqrt V(\Omega-x_1)$ and
\begin{equation*}
v(x)=\frac{1}{V}e_\e (u_\eps)(x_1+\frac{x}{\sqrt V})\quad\text{for }x\in B_{\mathcal{M}/2}\cap\widetilde\Omega, 
\end{equation*}
with the convention $B_r:=B_r(0)$ and we find 
\begin{gather}
1=v(0)\leq \sup_{B_{\mathcal{M}/2}\cap\widetilde\Omega} v \leq 4,\label{estimv1}\\
\frac1{s^{M-2}}\int_{B_s\cap\widetilde\Omega} v \leq E\leq\eta_0\quad\text{for all }s\leq \mathcal{M}/2.\label{estimv2}
\end{gather}
Since the Bochner-type inequality \eqref{est:Bochner} holds for $ u$ in $B_{r}(x_0) \cap \Omega$, we have
\begin{equation}\label{bochnerv}
-\Delta v \leq 8 v^2 \leq 32\, v \quad\text{in }B_{\mathcal{M}/2}\cap\widetilde\Omega.
\end{equation}
To prove \eqref{Eq:MGoal}, we consider separately the case $\mathcal{M} > 2$ and $\mathcal{M} \leq 2$.

\bigskip
\noindent
{\it Case 1:} $\mathcal{M} > 2$.
In this case, \eqref{estimv2} holds for every $s\in (0,1)$. If $V\leq 1$, then $\mathcal{M}^2 \leq 4s_1^2 \leq r^2\leq r_0^2$ for every $r\in (0, r_0)$ and then by choosing $r_0\leq \delta$,  \eqref{Eq:MGoal} follows.
It remains the case $V>1$. We show that $V$ is bounded by a constant $V_0 := C_0 \|\nabla_\tau u_{bd}\|_{C^{0,\beta}(\partial\Omega)}^2$ with $C_0$ depending only on $M$, $\beta$, $\kappa$ and $\Omega$ (defined in \eqref{Eq:eVd} below)  , which then implies $\mathcal{M}^2 = 4s_1^2 V \leq V_0r^2$, which proves \eqref{Eq:MGoal}.

For that, fix some arbitrary $p \in (M,\infty)$, e.g. $p = 2M$. Let $C_1 = C_1(M, \beta, p,\Omega)$ 
denote the constant $C$ in Lemma \ref{lem:bdryestim}. By \cite[Theorem 4.1]{HanLin}, there exists $C_2 = C_2(M)>0$  such that for any function $w \in H^1(B_1)$  satisfying $- \Delta w^+ \leq 32 w^+$ in $B_1=B_1(0)\subset \R^M$ then
\begin{equation}
w^+(0) \leq \frac{C_2}{|B_r|} \int_{B_r} w^+\,dx  \textrm{ for any } r \in (0,1). 	\label{Eq:DGMN}
\end{equation}

Note that if we introduce $\tilde \e=\sqrt{V} \e$ and $\tilde u_{\tilde \e}(x)=u_\e(x_1+V^{-1/2}x)$, then $-\Delta \tilde u_{\tilde \e}=\frac1{\tilde \e^2} W'(1-|\tilde u_{\tilde \e}|^2)\tilde u_{\tilde \e}$ in $\tilde \Omega$ and $v=e_{\tilde \e}(\tilde u_{\tilde \e})$. Thus, by applying Lemma~\ref{lem:bdryestim} to $\tilde \e$ and $\tilde u_{\tilde \e}$ and using \eqref{estimv1} and \eqref{estimv2}, we have\footnote{Note that we use that $V>1$ in the term ${V^{-1}}\|\nabla_\tau u_{bd}\|_{C^{0,\beta}(\partial\Omega)}^2$ as $V^{-1-\beta}<V^{-1}$.}
\begin{align}
\sup_{B_{1/2}\cap\partial\widetilde\Omega} v 
	& \leq C_1 \left( \norm{v}_{L^p(B_1\cap\widetilde\Omega)}^2 + {\norm{v}_{L^1(B_1\cap\widetilde\Omega)} + {V^{-1}}\|\nabla_\tau u_{bd}\|_{C^{0,\beta}(\partial\Omega)}^2    +\delta^2 }\right)\nonumber\\
	&\leq C_1 \left(4^{\frac{2(p-1)}{p}} \eta_0^{\frac{2}{p}} + \eta_0 + {V^{-1}}\|\nabla_\tau u_{bd}\|_{C^{0,\beta}(\partial\Omega)}^2    +\delta^2 \right).
		\label{Eq:Step2sup}
\end{align}
We choose some sufficiently small $\eta_0 , \delta > 0$ and sufficiently large $C_0 \geq 1$ (all depending only on $M$, $\beta$, $p$, $\kappa$, $\Omega$) such that 
\begin{equation}
C_1 \left(4^{\frac{2(p-1)}{p}} \eta_0^{\frac{2}{p}} + \eta_0 + {C_0^{-1}}    +\delta^2 \right) \leq \frac{1}{2(C_2 + 1)} < \frac{1}{2}.
	\label{Eq:eVd}
\end{equation}
Now, assume by contradiction that $V > V_0 = C_0 \|\nabla_\tau u_{bd}\|_{C^{0,\beta}(\partial\Omega)}^2$. Then by \eqref{Eq:Step2sup} and \eqref{Eq:eVd}, $v < 1/2$ on $B_{1/2}\cap\partial\widetilde\Omega$. By \eqref{bochnerv} and Kato's inequality (see \cite{Brezis84, Kato72}), 
the function
\begin{equation*}
\widetilde v:=\begin{cases}
\max\{v-1/2,0\}&\text{ in }B_{1/2}\cap\widetilde\Omega,\\
0 &\text{ in }B_{1/2}\setminus\widetilde\Omega,
\end{cases}
\end{equation*}
satisfies $-\Delta \widetilde v \leq 32\, \widetilde v$ in $B_{1/2}$ and the elliptic estimate \eqref{Eq:DGMN} yields (as $v(0)=1$)
\begin{equation}
1/2 = \widetilde v(0)\leq \frac{ C_2}{|B_{1/2}|} \int_{B_{1/2}}\widetilde v\leq \frac{ C_2}{|B_{1/2}|} \eta_0,
	\label{Eq:Step2Har}
\end{equation}
which implies a contradiction after possibly shrinking $\eta_0$. We deduce that $V \leq V_0$ and concludes the proof in Case 1.

\bigskip
\noindent
{\it Case 2:} $\mathcal{M} \leq 2$.
If $V \leq V_0=C_0 \|\nabla_\tau u_{bd}\|_{C^{0,\beta}(\partial\Omega)}^2$, then as in Case 1, we have $\mathcal{M}^2 = 4s_1^2 V \leq V_0r^2$ which yields \eqref{Eq:MGoal}. It remains to treat the case $V > V_0$. For that, applying Lemma~\ref{lem:bdryestim} for some $p\in (M, \infty)$ and using \eqref{estimv1} and \eqref{estimv2} give
\begin{align*}
\sup_{B_{\mathcal{M}/4}\cap\partial\widetilde\Omega} v 
	&\leq C_1\left((\mathcal{M}/2)^{2 - \frac{2M}{p}} \|v\|_{L^p(B_{\mathcal{M}/2} \cap \widetilde\Omega)}^2 
	+ (\mathcal{M}/2)^{-M} \|v\|_{L^1(B_{\mathcal{M}/2} \cap \widetilde\Omega)}
	+ {V^{-1}}\|\nabla_\tau u_{bd}\|_{C^{0,\beta}(\partial\Omega)}^2
	+ \frac{4\delta^2}{\mathcal{M}^2}
	\right)
	\\
	&\leq C_1\left(4^{\frac{2(p-1)}{p}} (\mathcal{M}/2)^{2 - \frac{4}{p}} \eta_0^{\frac{2}{p}} + {V^{-1}}\|\nabla_\tau u_{bd}\|_{C^{0,\beta}(\partial\Omega)}^2 + \frac{4(E + \delta^2)}{\mathcal{M}^2} \right).
\end{align*}
Since $\mathcal{M} \leq 2$, $2 - \frac{4}{p} > 0$ and $V > V_0$, we see from the above and \eqref{Eq:eVd} that
\begin{align*}
\sup_{B_{\mathcal{M}/4}\cap\partial\widetilde\Omega} v 
	&\leq \frac{1}{2(C_2 + 1)} +   \frac{4C_1(E + \delta^2)}{\mathcal{M}^2}  .
\end{align*}
For $|x| < \mathcal{M}/4$, define
\begin{equation*}
\widehat v(x) :=\begin{cases}
\max\Big\{v (x), \frac{1}{2(C_2 + 1)} +   \frac{4C_1(E + \delta^2)}{\mathcal{M}^2}\Big\}&\text{ in } x \in B_{\mathcal{M}/4}\cap\widetilde\Omega,\\
\frac{1}{2(C_2 + 1)} +   \frac{4C_1(E + \delta^2)}{\mathcal{M}^2} &\text{ in } B_{\mathcal{M}/4}\setminus\widetilde\Omega.
\end{cases}
\end{equation*}
By \eqref{bochnerv} and Kato's inequality, $-\Delta \widehat v \leq   32  \widehat v$ in $B_{\mathcal{M}/4}$, and the elliptic estimate \eqref{Eq:DGMN} gives
\begin{equation*}
1 = v(0)  \leq \widehat v(0)\leq \frac{ C_2}{|B_{\mathcal{M}/4}|} \int_{B_{\mathcal{M}/4}}\widehat v \leq \frac{ C_2}{|B_{\mathcal{M}/4}|} \int_{B_{\mathcal{M}/4} \cap \widetilde \Omega}  v + \underbrace{\frac{C_2}{2(C_2 + 1)}}_{\leq 1/2} +   \frac{4C_2 C_1(E + \delta^2)}{\mathcal{M}^2}.
\end{equation*}
Recalling \eqref{estimv2}, we thus obtain
\[
\frac{1}{2} \leq  \frac{4C_2 (C_1 + 4|B_1|^{-1})(E + \delta^2)}{\mathcal{M}^2},
\]
that is
\[
\mathcal{M}^2   \leq8C_2(C_1 + 4|B_1|^{-1})(E + \delta^2)  .
\]
We have thus established \eqref{Eq:MGoal} in this case too. The proof is complete.
\end{proof}

\bigskip

Now we are in position to prove the uniform gradient estimate in Proposition \ref{pro:main} for our solutions $u_\eps$ in the region where their 
limit $u_0$ is continuous.

\begin{proof}[Proof of Proposition \ref{pro:main}]
By Remark \ref{rem:regu} and Theorem \ref{thm:MPI}, we know that $u_\eps\in C^1(\bar \Omega)$ and $|u_\eps|\leq 1$ in $\Omega$. Let $r_0$, $\delta$ and $\eta_0$ be as in Proposition \ref{prop:smallestim}. By hypothesis, after possibly shrinking $r_0$, we may assume that $u_0$ is continuous in
\[
K' = \bar\Omega \cap (\bar B_{r_0}(0)+K) =\{x+y \in \bar \Omega\, :\, |x|\leq r_0, y\in K\}.
\]
Moreover, by Lemma \ref{lemma:unifconvmodule}, we can find $\eps_0 > 0$ such that, after possibly shrinking $r_0$ further,
\[
1 - |u_\eps| \leq \delta \text{ in } K' \text{ for } 0 < \eps \leq \eps_0.
\]
By Lemma \ref{lemma:sEc}, there exists $0 < r_1 \leq r_0$ such that, after possibly shrinking $\eps_0$, it holds that
\[
\frac{1}{s^{M-2}} \int_{\Omega \cap B_s(x)} e_\eps(u_\eps) \leq \eta_0 \quad \forall~x \in K, s \in (0,r_1], \eps \in (0,\eps_0].
\]
Applying Proposition \ref{prop:smallestim}, we reach the desired gradient estimate for $\eps \leq \eps_0$. When $\eps > \eps_0$, the desired gradient estimate follows from \eqref{eq:bulk3}.
\end{proof}

\section{Uniform Laplacian estimate}\label{Sec:LapEst}

The main result of this section is the uniform Laplacian estimate for solutions $u_\eps$ to \eqref{eq:eps}
converging as in \eqref{ass:convH1} to a stationary harmonic map $u_0$ in the region where $u_0$ is continuous.

\begin{proposition} 
\label{pro:main1}
Let $W \in C^2((-\infty,1])$ be such that \eqref{ass:Wstructure}, \eqref{ass:Wnondeg} and \eqref{ass:convex} hold, $\Omega\subset \R^M$ be a $C^2$ bounded domain and 
$u_{bd}\in C^{1, \beta}(\partial \Omega, \Ss^{N-1})$ for some $\beta\in (0,1)$.
Let $(u_\eps)_{\eps\to 0}\subset H^1\cap L^\infty(\Omega, \R^N)$ be a family of solutions to \eqref{eq:eps}  that  satisfies \eqref{ass:convH1} for a stationary harmonic map $u_0\in H^1(\Omega, \Ss^{N-1})$. Then for any compact $K\subset\overline\Omega$ where $u_0$ is continuous in a neighborhood of $K$, there exists a constant $C > 0$ depending only on $M, N, W, \Omega, K, u_{bd}$ and $u_0$ such that
$$
\sup_{K}|\Delta u_{\eps}|\le C \quad \forall~\eps >0.$$
Moreover, 
\[
\|u_\vareps\|_{C^{1,\beta}(K)} \leq C \quad \forall~\eps >0.
\]
\end{proposition}

Clearly, Proposition \ref{pro:main1} implies Theorem \ref{thm:MDI} in the introduction.
 
The ideas of proof of Proposition \ref{pro:main1} are as follows: Using the non-degeneracy assumption \eqref{ass:Wnondeg} and the convexity assumption \eqref{ass:convex}, it is shown that the function $W(1 - |u_\eps|^2)$ satisfies a partial differential inequality (see \eqref{PDE-inequa}) of the form:
\begin{equation}
-\eps^2 \Delta f + a f^{\frac{2\alpha}{\alpha + 1}} \leq C,
	\label{Eq:fDI}
\end{equation}
with $a > 0, C> 0$.
When $\alpha = 1$, \eqref{Eq:fDI} is linear and one can use arguments from \cite{BBH93, NZ13} to obtain a pointwise upper bound for $W(1 - |u_\eps|^2)$. When $\alpha > 1$, \eqref{Eq:fDI} is semi-linear and we borrow ideas available in the context of conformal geometry \cite{LoewnerNirenberg} to prove a pointwise upper bound for $W(1 - |u_\eps|^2)$. Using the non-degeneracy assumption \eqref{ass:Wnondeg} and the bound for $W(1 - |u_\eps|^2)$, we deduce a bound for $W'(1 - |u_\eps|^2)$ and hence a bound for $\Delta u_\eps$ (by \eqref{eq:eps}), as desired.

We start with $L^\infty$ estimates for \eqref{Eq:fDI} in Lemmas \ref{BBHMaxPrin::Bdry} and \ref{BBHNLMaxPrin::Bdry}. In the linear case, we will use the fact that, for $\gamma \geq 0$, the modified Bessel equation
$$
r^2y''+r y'-(r^2+\gamma^2)y=0 \quad \textrm{ in } (0,\infty).
$$
has a unique, up to a multiplicative scalar, solution which is bounded as $r \rightarrow 0$, namely the modified Bessel function of the first kind $I_\gamma:(0, \infty)\to (0, \infty)$. Moreover, $I_\gamma(r) \sim r^\gamma$ as $r \rightarrow 0$ and $I_\gamma(r) \sim e^r\,r^{-1/2}$ as $r \rightarrow \infty$. See e.g. \cite[\S 9.6.1, \S 9.6.7, \S 9.6.10 and \S 9.7.1]{AbSt}. In particular, there exists a constant $d_\gamma \geq 1$ such that
\begin{equation}
\frac{1}{d_\gamma} \leq \frac{I_\gamma(r/2) e^{r/2}}{I_\gamma(r)} \leq d_\gamma \text{ for all } r \in (0,\infty).
	\label{Eq:BIdoubling}
\end{equation}

\begin{lemma}\label{BBHMaxPrin::Bdry}
Let $M\ge 2$, $\Omega$ be a domain in $\RR^M$, $x_0$ $\in$ $\bar\Omega$ and $R > 0$ such that $B_R(x_0) \cap \Omega$ is a Lipschitz domain. Assume that $f \in H^1(B_R(x_0) \cap \Omega)$ satisfies in the weak sense
$$-\eps^2\,\Delta f + a\,f \leq C \text{ in } B_R(x_0) \cap \Omega,
$$
where $\eps$, $a$ and $C$ are positive constants. If $B_R(x_0) \cap \partial \Omega$ is non-empty, assume also that
\[
f = 0 \text{ on } B_R(x_0) \cap \partial\Omega.
\]
With 
$\gamma=\frac{M-2}{2}\geq 0$, we have
\[
f(x) \leq \frac{C}{a} + \Big[\sup_{\partial B_R(x_0) \cap \Omega}\Big(f - \frac{C}{a}\Big)^+\Big]\,\frac{I_\gamma\left(\frac{\sqrt{a}|x-x_0|}{\eps}\right) R^{\gamma}}{I_\gamma\left(\frac{\sqrt{a}R}{\eps}\right)|x-x_0|^{\gamma}}
\quad \textrm{ for almost all  }x \in B_R(x_0) \cap \Omega.
\] 
Consequently, 
\[
f(x) \leq \frac{C}{a} + 2^{\gamma}d_\gamma\Big[\sup_{\partial B_R(x_0) \cap \Omega}\Big(f - \frac{C}{a}\Big)^+\Big] \exp\Big(-\frac{\sqrt{a}R}{2\eps}\Big) \quad \text{ for almost all } x\in B_{R/2}(x_0) \cap \Omega.
\]
\end{lemma}

\begin{proof} 
Without loss of generality, we assume $x_0 = 0$ and write $B_R$ in place of $B_R(x_0)$.

 We first find a barrier function $\phi$ in $B_R$   that will be used for comparison with $f$. For that, we proceed as in the proof of \cite[Lemma $6$]{NZ13} by imposing that  $\phi|_{\partial B_R}\equiv 1$, $\phi\geq 0$ in $B_R$  and 
 \be
 \label{PDE-phi}
 -\eps^2\,\Delta\phi + a\,\phi=0 \quad \text{ in } B_R.
 \ee
We look for such $\phi$ in the class of radial functions, i.e., $\phi(x)=\phi(r)$, $r=|x|\in (0, R)$ which means that $\phi$  solves the ODE 
$$
-\eps^2(\phi''+\frac{M-1}{r}\phi')+a\phi=0 \quad \textrm{in } (0,R).
$$ 
The change of variables $\phi(r):=r^{-\gamma} y(r)$ with $\gamma=\frac{M-2}{2}\geq 0$ yields a modified Bessel equation for the function $y$:
$$r^2y''+ry'-(\frac{a}{\eps^2}r^2+\gamma^2)y=0 \quad \textrm{in } (0,R).$$ We are thus led to set
\[
y(r) = I_\gamma(\frac{\sqrt{a}}{\eps}r)
\]
and
\be
\label{def-phi}
\phi(x):=\frac{I_\gamma\left(\frac{\sqrt{a}|x|}{\eps}\right) R^{\gamma}}{I_\gamma\left(\frac{\sqrt{a}R}{\eps}\right)|x|^{\gamma}} >0 \quad \text{ for every } x\in B_R\setminus\{0\}.
\ee
As $I_\gamma(s) \sim s^\gamma$ as $s \rightarrow 0$, we have $\phi\in L^\infty(B_R )$. Since $\phi$ satisfies the PDE \eqref{PDE-phi} in $B_R \setminus \{0\}$, standard elliptic regularity yields that $\phi$ $\in$ $C^\infty(B_R)$ and $\phi$ satisfies \eqref{PDE-phi} in all of $B_R$. 

Now the function
\[
\tilde f(x) := f(x) -  \frac{C}{a}  - \Big[\sup_{ \partial B_R \cap \Omega }\Big(f - \frac{C}{a}\Big)^+\Big]\,\phi(x) 
\]
satisfies $-\eps^2 \Delta \tilde f + a\tilde f \leq 0$ in $B_R \cap \Omega$ and $\tilde f \leq  0$ on $\partial (B_R \cap \Omega)$. By the weak maximum principle, $\tilde f \leq 0$ in $B_R \cap \Omega$, which yields the first conclusion.

Next, note that since $(r^{M-1}\phi')' = \frac{a}{\eps^2} r^{M-1} \phi \geq 0$ (by \eqref{PDE-phi}) and $\phi'(0) = 0$ (by smoothness and symmetry of $\phi$ in $B_R$) we have $\phi' \geq 0$. Hence 
\[
\phi(x) \leq \phi(R/2) = \frac{2^{\gamma}I_\gamma\left(\frac{\sqrt{a}R}{2\eps}\right) }{I_\gamma\left(\frac{\sqrt{a}R}{\eps}\right) } \text{ for }x \in B_{R/2}.
\]
Using \eqref{Eq:BIdoubling} we then have
\[
\phi(x) \leq \phi(R/2) \leq  2^{\gamma}d_\gamma\exp\Big(-\frac{\sqrt{a}R}{2\eps}\Big) \text{ for } x \in B_{R/2}.
\]
The second conclusion follows from the above estimate and the first conclusion.
\end{proof}

We next turn to the semi-linear case.

\begin{lemma}\label{BBHNLMaxPrin::Bdry}
Let $M\ge 2$, $\gamma > 1$, $\Omega$ be a domain in $\RR^M$, $x_0$ $\in$ $\bar\Omega$, and $R > 0$ such that $B_R(x_0) \cap \Omega$ is a Lipschitz domain. Assume that $f \in H^1 \cap L^\infty(B_R(x_0) \cap \Omega) $ is non-negative and satisfies in the weak sense
$$-\eps^2\,\Delta f + a\,f^{\gamma} \leq C \text{ in } B_R(x_0) \cap \Omega,
$$
where $\eps$, $a$ and $C$ are positive constants. If $B_R(x_0) \cap \partial \Omega$ is non-empty, assume also that
\[
f = 0 \text{ on } B_R(x_0) \cap \partial\Omega.
\]
Then there exists $c = c(M,\gamma) > 0$ such that
\[
f(x) \leq \Big(\frac{C}{a}\Big)^{\frac{1}{\gamma}} + \Big(\frac{cR^2 \eps^2}{a}\Big)^{\frac{1}{\gamma - 1}}(R^2 - |x - x_0|^2)^{-\frac{2 }{\gamma-1}} \quad \textrm{ for }x \in B_R(x_0) \cap \Omega.
\] 
\end{lemma}

\begin{proof} 
Without loss of generality, we assume $x_0 = 0$ and write $B_R$ in place of $B_R(x_0)$. Replacing $\Omega$ by $\Omega \cup B_R$ and extend $f$ by zero in $B_R \setminus \Omega$, we may assume without loss of generality that $B_R \subset \Omega$. (Here we have used that $f$ is non-negative.)

We start by constructing positive barrier functions $Z$ in $B_R$:
\begin{equation}
\begin{cases}
-\eps^2\,\Delta Z + a\,Z^{ \gamma} \geq C \text{ in } B_R,\\
Z(x) \rightarrow \infty \text{ as } |x|  \rightarrow R.
\end{cases}
	\label{Eq:ZNLB}
\end{equation}
This resembles the Loewner-Nirenberg equation in conformal geometry \cite{LoewnerNirenberg}. Motivating by this work, we let $r = |x|$ and use the ansatz
 $$Z(x) = Z_{\lambda,\mu}(x):=\lambda( R^2-r^2)^{-\frac{2 }{\gamma-1}}+\mu \quad \textrm{ for } x\in B_R$$ with $\lambda$ and $\mu$ positive, to be fixed later. Since $\lambda > 0$, $Z(x) \rightarrow \infty$ as $ |x|  \rightarrow R$. We compute
\begin{align*} 
\Delta Z&=\frac{1}{r^{M-1}} \left( r^{M-1} Z'(r) \right)' \\
&= \frac{4\lambda }{(\gamma-1)r^{M-1}}\left(r^M(R^2-r^2)^{-\frac{\gamma + 1}{\gamma-1} }\right)' \\
             &= \frac{4\lambda }{ \gamma-1 } (R^2-r^2)^{-\frac{2\gamma }{\gamma-1} }\left[M(R^2-r^2)+\frac{2(\gamma+1)}{\gamma-1}r^2\right]\\
             &= \frac{4\lambda }{ \gamma-1 } (R^2-r^2)^{-\frac{2\gamma }{\gamma-1}}\left[MR^2-(M-\frac{2(\gamma+1)}{\gamma-1})r^2\right].
\end{align*} 
Using the inequality $(s_1+s_2)^\gamma\geq s_1^\gamma+s_2^\gamma$ for  $s_1, s_2\geq 0$, we also have
\begin{align*} 
Z^{\gamma}
	= \left[\lambda(R^2-r^2)^{-\frac{2 }{\gamma-1}}+\mu\right]^{\gamma}
	\ge \lambda^{\gamma}(R^2-r^2)^{-\frac{2\gamma}{\gamma-1}}+\mu^{\gamma}.
\end{align*} 
 It follows that
\begin{align*}
&-\eps^2\Delta Z+ aZ^{\gamma}\\
	&\qquad \ge - \frac{4\vareps^2\lambda }{ \gamma-1 } (R^2-r^2)^{-\frac{2\gamma }{\gamma-1}}\left[MR^2-(M-\frac{2(\gamma+1)}{\gamma-1})r^2\right]
		+ a \left[ \lambda^{\gamma}(R^2-r^2)^{-\frac{2\gamma}{\gamma-1}}+\mu^{\gamma}\right]\\
&\qquad \geq \frac{4\vareps^2\lambda }{ \gamma-1 } (R^2-r^2)^{-\frac{2\gamma }{\gamma-1}}  \underbrace{\left[\frac{a(\gamma-1)}{4\eps^2}\lambda^{\gamma-1}
	- MR^2 + (M-\frac{2(\gamma+1)}{\gamma-1})r^2\right]}_{=:I}
+ a\mu^{\gamma}.
\end{align*}
In order for $Z$ to satisfy \eqref{Eq:ZNLB}, the right hand side above needs to be greater than $C$. To this end, we pick 
\[
\lambda = \Big[\frac{4 }{\gamma - 1}\Big(M + \frac{2(\gamma + 1)}{\gamma - 1}\Big)\frac{R^2 \eps^2}{a}\Big]^{\frac{1}{\gamma - 1}} \quad (\text{so that }I > 0),
\]
and impose that
\[
a\mu^{\gamma} \geq C \quad \Leftrightarrow \quad \mu \geq \Big(\frac{C}{a}\Big)^{\frac{1}{\gamma}} =: \underline{\mu}.
\] 

With $\lambda$ and $\underline{\mu}$ as above, to finish the proof, we only need to show that
\begin{equation}
f \leq Z_{\lambda,\underline{\mu}} \text{ a.e. in } B_R.
	\label{Eq:fZbClaim}
\end{equation}
Suppose by contradiction that
\[
z := \inf_{B_R} (Z_{\lambda,\underline{\mu}} - f) < 0.
\]
(Note that since $f \in L^\infty(B_R)$, $z$ is finite.) Let $\mu = \underline{\mu} - z > \underline{\mu}$ and $g = Z_{\lambda,\mu} - f$. Then $g \in H^1_{\rm loc} \cap L^\infty_{\rm loc}(B_R)$, $\inf_{B_R} g = 0$ and $g$ satisfies
\[
-\eps^2\Delta g+a(x) g \ge 0 \text{ in } B_R,
\]
where
$$a(x)=\begin{cases}
  \frac{a(Z_{\lambda,\mu }(x)^{\gamma}-f(x)^{\gamma})}{Z_{\lambda,\mu }(x)- f(x)} & \textrm{ if } Z_{\lambda,\mu }(x) > f(x)\\
a \gamma Z_{\lambda,\mu }(x)^{\gamma-1} & \textrm{ if } Z_{\lambda,\mu }(x)=f(x).
 \end{cases}
 $$
Since $Z_{\lambda,\mu}(x) \rightarrow \infty$ as $|x| \rightarrow R$ and $f \in L^\infty(B)$, there exists $0 < R' < R$ such that
\[
\inf_{B_R \setminus B_{R'}} g \geq  1 \text{ and hence }
\inf_{B_{R'}} g = \inf_{B_R} g = 0.
\]
By the weak Harnack inequality (see e.g. \cite[Theorem 8.18]{GT}), this implies that  
\[
\|g\|_{L^1(B_{(R+R')/2})} \leq C(M, a, \lambda, \mu, R, R') \inf_{B_{R'}} g = 0,
\]
i.e. $g = 0$ a.e. in $B_{(R+R')/2}$. It follows that
\[
0 =  \|g\|_{L^\infty(B_{(R+R')/2} )} \geq \inf_{B_{(R+R')/2} \setminus B_{R'}} g   \geq \inf_{B_R \setminus B_{R'}} g \geq   1,
\]
which is absurd. We have thus proved \eqref{Eq:fZbClaim} and concluded the proof.
\end{proof}

\bigskip

We now prove a localised version of the uniform Laplacian estimate claimed in Proposition \ref{pro:main1}:
\begin{lemma}\label{bd:unifLap}
Let $W \in C^2((-\infty,1])$ be such that \eqref{ass:Wstructure}, \eqref{ass:Wnondeg} and \eqref{ass:convex} hold, $\Omega\subset \R^M$ be a $C^2$ bounded domain and 
$u_{bd}\in C^{1, \beta}(\partial \Omega, \Ss^{N-1})$ for some $\beta\in (0,1)$. Then there exist $r_1, \delta_1, C_1 >0$ (depending only on $M$, $N$, $W$ and $\Omega$) such that 
for every $\eps > 0$,  $r\in (0,r_1)$, $x_0\in\overline\Omega$ and any solution $u_\eps\in H^1\cap L^\infty(\Omega, \R^N)$ to \eqref{eq:eps}  satisfying
$$\textrm{$1-|u_\eps| <\delta_1$  in $B_r(x_0) \cap \Omega$},$$
we have
\[
\|\Delta u_\eps\|_{L^\infty(B_{r/2}(x_0) \cap \Omega)} \leq C_1\big(\|\nabla u_\eps\|_{L^\infty(B_r(x_0) \cap \Omega)}^2 + r^{-2}\big) .
\]
\end{lemma}

\bproof To simplify the presentation, we drop the index $\eps$. We denote
$$
Y(x):=\frac{W(1-|u(x)|^2)}{2}. 
$$ 
Recall that $u\in C^2(\Omega)\cap  C^{1,\beta}(\Omega)$ (see Remark \ref{rem:regu}) and $|u|\leq 1$ in $\Omega$ by Theorem \ref{thm:MPI}. 

We start with the choice of $r_1$ and $\delta_1$.  We fix $r_1 > 0$ sufficiently small such that $B_s(x) \cap \Omega$ is a Lipschitz domain for all $s \in (0,r_1)$ and $x \in \bar\Omega$ with $\mathrm{dist}(x,\partial\Omega) \leq 2s/3$. (See the set up in the proof of Corollary \ref{Cor:BRE}.) This choice of $r_1$ will become clear later on when we apply Lemmas \ref{BBHMaxPrin::Bdry} and \ref{BBHNLMaxPrin::Bdry}.

By \eqref{ass:Wnondeg}, there exists $0 < \kappa' \leq \kappa$ such that
\[
\begin{cases}
\frac{1}{2} c_0(\alpha + 1) t^\alpha \leq W'(t) \leq 2c_0(\alpha + 1) t^\alpha,\\
\frac{1}{2}c_0 t^{\alpha + 1} \leq W(t) \leq 2c_0 t^{\alpha + 1}
\end{cases} \text{ for } t \in [0, \kappa'].
\]
In particular,
\begin{equation}
2^{-\frac{2\alpha+1}{\alpha+1}}c_0^{\frac{1}{\alpha+1}} (\alpha + 1)W(t)^{\frac{\alpha}{\alpha+1}} \leq W'(t) \leq  2^{\frac{2\alpha+1}{\alpha+1}}c_0^{\frac{1}{\alpha+1}} (\alpha + 1)W(t)^{\frac{\alpha}{\alpha+1}} \text{ for } t \in [0,\kappa'].
	\label{Eq:WndC}
\end{equation}
We can take $\delta_1 = \min\{1/4, \kappa'/2\}$. 

We now fix a ball $B_r(x_0)$ with $x_0 \in \bar \Omega$ and $r \in (0,r_1)$ such that $1-|u_\eps| <\delta_1$  in $B_r(x_0) \cap \Omega$ and proceed to bound $\Delta u_\eps$ in $B_{r/2}(x_0) \cap \Omega$. In the proof,  $C$ denotes a generic positive constant which depends only on the dimension $M$, $N$, $r_1$, $\delta_1$, and the behaviour of the function $W$ in $[0, 1]$. 

Note that, by \eqref{eq:eps},
\[
|\Delta u_\eps| \leq  \frac{1}{\eps^2} \max_{[0,1]} |W'| \text{ in }  \Omega.
\]
Hence if $r \leq \eps$, then the conclusion follows. We assume in the rest of the proof that $r > \eps$.

\medskip
\noindent{\it Case 1:} $\mathrm{dist}(x_0, \partial \Omega) \geq r$ or $\mathrm{dist}(x_0, \partial \Omega) \leq 2r/3$. In particular, $B_r(x_0) \cap \Omega$ is a Lipschitz domain. Without loss of generality, we can take $x_0 = 0$ and write $B_r$ in place of $B_r(x_0)$.

We compute 
\bea
\nonumber
\Delta Y&=- \sum_i \nabla_i [W'(1-|u|^2) u \cdot (\nabla_i u)]\\
&= 2W''(1-|u|^2) u^2 |\nabla u|^2 - W'(1-|u|^2)|\nabla u|^2-
W'(1-|u|^2)u \cdot (\Delta u)\nonumber\\
&= 2W''(1-|u|^2) u^2 |\nabla u|^2 - W'(1-|u|^2)|\nabla u|^2 + \frac{1}{\eps^2}
(W'(1-|u|^2))^2 |u|^2,
\label{est:LapY}
 \eea
 where we have used \eqref{eq:eps} for the last identity.
 
 Note that $0 \leq 1 -|u|^2 = (1 + |u|) (1 - |u|) \leq 2\delta_1 \leq  \kappa' \leq \kappa$ in $B_r \cap \Omega$. Using this as well as \eqref{ass:convex} and \eqref{Eq:WndC}, we can estimates the terms on the right hand side of \eqref{est:LapY} as follows. First, by the convexity assumption \eqref{ass:convex}, we have
  $$
 W''(1-|u|^2) u^2 |\nabla u|^2 \ge 0 \textrm{ in } B_r  \cap \Omega.
  $$ 
Second, using the second inequality in \eqref{Eq:WndC}, we have
\[
-W'(1-|u|^2)|\nabla u|^2 \ge - C \|\nabla u\|_{L^\infty(B_r  \cap \Omega)}^2 Y^{\frac{\alpha}{\alpha+1}} \textrm{ in } B_r  \cap \Omega \text{ and for } 1 \leq i \leq M.
 \]
 Third, using $|u| \leq 1$ and the first inequality in \eqref{Eq:WndC}, we have
 $$
 (W'(1-|u|^2))^2 |u|^2 \ge \frac{1}{C} Y^{\frac{2\alpha}{\alpha+1}} \textrm{ in $B_r \cap \Omega$}.
 $$
Putting the above estimates into \eqref{est:LapY}, we get for some $\hat{C}> 1$ depending only on $W$:
 \[
 \eps^2\Delta Y\ge \frac{1}{C} Y^{\frac{2\alpha}{\alpha+1}}- \eps^2 C \|\nabla u \|_{L^\infty(B_r  \cap \Omega)}^2 Y^{\frac{\alpha}{\alpha+1}}
\ge \frac{1}{\hat C}Y^{\frac{2\alpha}{\alpha+1}}-\eps^4\hat{C}\|\nabla u \|_{L^\infty(B_r  \cap \Omega)}^4 \textrm{ in } B_r \cap \Omega.
 \]
 Equivalently, we have
  \be
 \label{PDE-inequa}
- \eps^2\Delta Y + \frac{1}{\hat C}Y^{\frac{2\alpha}{\alpha+1}} \leq\eps^4\hat{C} \|\nabla u \|_{L^\infty(B_r  \cap \Omega)}^4 \textrm{ in } B_r \cap \Omega.
 \ee

\nd {\it Subcase 1a:} $\alpha=1$. Using \eqref{PDE-inequa} and applying  Lemma ~\ref{BBHMaxPrin::Bdry}, we obtain
\[
0 \leq Y \leq \eps^4 \hat C^2 \|\nabla u \|_{L^\infty(B_r  \cap \Omega)}^4 + C \exp\Big(-\frac{r}{2\sqrt{\hat C}\eps}\Big) \leq C \eps^4( \|\nabla u \|_{L^\infty(B_r  \cap \Omega)}^4 + r^{-4}) \textrm{ in } B_{r/2} \cap \Omega.
\]
Using \eqref{eq:eps} and the fact that $r > 0$ and noting that $|W'(1-|u|^2)| \leq CY^{\frac{\alpha}{\alpha+1}} = C Y^{1/2}$ (by \eqref{Eq:WndC}), we deduce the desired estimate $\|\Delta u\|_{L^\infty(B_{r/2} \cap \Omega)} \leq C ( \|\nabla u \|_{L^\infty(B_r  \cap \Omega)}^2 + r^{-2})$.

 \bigskip
 
\nd {\it Subcase 1b:} $\alpha>1$. Using \eqref{PDE-inequa} and applying  Lemma ~\ref{BBHNLMaxPrin::Bdry} (with $\gamma = \frac{2\alpha}{\alpha + 1} > 1$), we obtain
\begin{align*}
0 \leq Y 
	&\leq \Big(\eps^4 \hat C^2 \|\nabla u \|_{L^\infty(B_r  \cap \Omega)}^4\Big)^{\frac{\alpha + 1}{2\alpha}} + \Big(\frac{C \hat C  \eps^2}{r^2}\Big)^{\frac{\alpha+1}{\alpha - 1}} \\
	&\leq C \eps^{\frac{2(\alpha+1)}{\alpha}} \Big( \|\nabla u \|_{L^\infty(B_r  \cap \Omega)}^{\frac{2(\alpha + 1)}{\alpha}} + \eps^{\frac{2(\alpha+1)}{\alpha(\alpha - 1)}}r^{-\frac{2(\alpha + 1)}{\alpha - 1}}\Big)\quad \textrm{ in } B_{r/2} \cap \Omega.
\end{align*}
In particular, since $r > \eps$, we have
\[
0 \leq Y \leq C \eps^{\frac{2(\alpha+1)}{\alpha}} ( \|\nabla u \|_{L^\infty(B_r  \cap \Omega)}^{\frac{2(\alpha + 1)}{\alpha}} + r^{-\frac{2(\alpha + 1)}{\alpha}})\quad \textrm{ in } B_{r/2} \cap \Omega.
\]
Recalling \eqref{eq:eps}, noting that $|W'(1-|u|^2)| \leq CY^{\frac{\alpha}{\alpha+1}}$ (by \eqref{Eq:WndC}) and using the inequaity $(a+b)^{\frac{\alpha}{\alpha+1}} \leq a^{\frac{\alpha}{\alpha+1}} +b^{\frac{\alpha}{\alpha+1}}$ for $a, b \geq 0$, we deduce the desired estimate $\|\Delta u\|_{L^\infty(B_{r/2} \cap \Omega)} \leq C ( \|\nabla u \|_{L^\infty(B_r  \cap \Omega)}^2 + r^{-2})$.

\bigskip
\noindent{\it Case 2:} $2r/3 < \mathrm{dist}(x_0, \partial \Omega) < r$. In this case $B_{r/2}(x_0) \subset \Omega$ and for every $x \in B_{r/2}(x_0)$, $B_{r/6}(x) \subset \Omega$. By Case 1, we have 
\[
|\Delta u(x)| \leq C ( \|\nabla u \|_{L^\infty(B_{r/6}(x))}^2 + r^{-2}) \leq C ( \|\nabla u \|_{L^\infty(B_r(x_0)  \cap \Omega)}^2 + r^{-2}) \text{ for all } x \in B_{r/2}(x_0).
\]
We again obtain $\|\Delta u\|_{L^\infty(B_{r/2}(x_0) \cap \Omega)} \leq C ( \|\nabla u \|_{L^\infty(B_r(x_0)  \cap \Omega)}^2 + r^{-2})$.
 \eproof

Now we prove the main result in Proposition \ref{pro:main1}.

\begin{proof}[Proof of Proposition \ref{pro:main1}]
By Remark \ref{rem:regu} and Theorem \ref{thm:MPI}, we know that $u_\eps\in C^2 (\Omega) \cap C^{1,\beta}(  \Omega)$ and $|u_\eps|\leq 1$ in $\Omega$. Let $r_1, \delta_1$ be as in Lemma \ref{bd:unifLap}. Let $\omega_1 \subset \omega_2 \subset \bar \Omega$ be neighborhoods of $K$ where $u_0$ is continuous such that, after possibly shinking $r_1$ slightly, $B_{r_1}(x)\cap \Omega\subset \omega_1\cap \Omega$ for every $x\in K$ and $B_{r_1}(x)\cap \Omega\subset \omega_2\cap \Omega$ for every $x\in \omega_1$. By Lemmas \ref{lemma:sEc} and \ref{lemma:unifconvmodule} and the continuity of $u_0$ on $\omega_2$, we know that $|u_\eps|\to 1$ uniformly in $\omega_2$. Hence there exists $\eps_0 >0$ so that
\[
0 \leq 1 - |u_\eps| \leq \delta_1 \text{ in } \omega_2 \text{ for all } \eps \in (0,\eps_0).
\]
In the sequel, the constant $C$ denotes a generic constant depending only on $M, N, W, \Omega, K$, $r_1, \delta_1, \eps_0, \omega_1, \omega_2, u_{bd}$ and $u_0$.

By Proposition \ref{pro:main}, after possibly shrinking $\eps_0$, we have
\[
|\nabla u_\eps|\leq  C \text{ in } \omega_2 \text{ for all } \eps \in (0,\eps_0). 
\]
Applying Lemma \ref{bd:unifLap} to $u_\eps$ in $B_{r_1}(x) \cap \Omega$ for any $x\in \omega_1$, we have
\[
|\Delta u_\eps|\leq  C \text{ in } \omega_1 \text{ for all } \eps \in (0,\eps_0). 
\]
Enlarging $C$ if necessary and recalling \eqref{eq:eps}, we see that this estimate also hold for $\eps \geq \eps_0$. Applying $W^{2,p}$-interior estimate for any $p < \infty$ (e.g. \cite{GT}[Theorem 9.11]) and Morrey's inequality, we obtain interior H\"older $C^{0, \beta}$ estimate for $\nabla u_\vareps$:
\[
\|\nabla u_\eps\|_{C^{0,\beta}(B_{r_1/4}(x_0))} \leq C \quad \forall x_0 \in K \text{ such that } B_{r_1/2}(x_0) \subset \omega_1.
\]
Applying Lemma \ref{lem:bdry-est}, we obtain boundary H\"older $C^{0, \beta}$ estimate for $\nabla u_\vareps$:
\[
\|\nabla u_\eps\|_{C^{0,\beta}(B_{r_1/2}(x_0) \cap \Omega)} \leq C \quad \forall x_0 \in K \cap \partial\Omega.
\]
These two estimates imply that
\[
\|\nabla u_\eps\|_{C^{0,\beta}(K)} \leq C.
\]
The proof is complete.
\end{proof}

\medskip
\noindent{\bf Acknowledgment.}  R.I.  is partially supported by the ANR projects ANR-21-CE40-0004 and
ANR-22-CE40-0006-01. A.Z. has been partially supported by the Basque Government through the BERC 2022-2025 program and by the Spanish State Research Agency through BCAM Severo Ochoa CEX2021-001142-S and through project PID2023-146764NB-I00 funded by MICIU/AEI/10.13039/501100011033 and cofunded by the European Union. A.Z. was also partially supported by a grant of the Ministry
of Research, Innovation and Digitization, CNCS-UEFISCDI, project number PN-IV-P2-2.1-T-TE-2023-
1704, within PNCDI IV.

\smallskip
\noindent{\bf Right retention statement. } For the purpose of open access, the authors
have applied a CC-BY public copyright licence to any author accepted
manuscript arising from this submission.

\appendix 

\section{Appendix}

\subsection{Weak and stationary harmonicity of the limit}\label{App:HM}

\begin{lemma}
\label{lem:Stationarity} Let $\Omega\subset \R^M$ be a domain, $W\in C^1((-\infty, 1])$ and $u_\eps\in H^1 \cap L^\infty(\Omega, \R^N)$  be a solution to \eqref{eq:eps} for small $\eps>0$. If $u_\eps\rightharpoonup u_0$ weakly in $H^1$ to a limit $u_0\in H^1(\Omega, \Ss^{N-1})$ as $\eps\to 0$, then $u_0$ is a weakly $\Ss^{N-1}$-valued harmonic map, i.e, $-\Delta u_0 =u_0 |\nabla u_0 |^2$ in the distribution sense in $\Omega$. If, in addition, we assume  that \eqref{ass:convH1} holds and $W$ satisfies \eqref{ass:Wsweak}, then $u_\eps \rightarrow u_0$ in $H^1(\Omega)$, $\frac{1}{\vareps^2}W(1 - |u_\vareps|^2) \rightarrow 0$ in $L^1(\Omega)$ and $u_0$ is a stationary harmonic map, i.e, the stress-energy tensor
\[
S_{ij} = \frac{1}{2} |\nabla u_0|^2 \delta_{ij} - \nabla_i u_0 \cdot \nabla_j u_0, \quad 1\leq i, j\leq M
\]
is divergence-free in the distributional sense in $\Omega$.
\end{lemma}

\begin{proof} To simplify the notation, in the following, we write $u$ instead of $u_0$.
 It is known that $u\in H^1(\Omega, \Ss^{N-1})$ is a weakly harmonic map in $\Omega$ if and only if  (see e.g.\cite[Lemma 3.7]{Mosbook}):
$$
\nabla \cdot  (u\wedge \nabla u)=0 \textrm{ in the distribution sense in } \Omega,
$$ where $\cdot$ is the scalar product in $\R^M$ and $\wedge$ denotes the wedge product in $\R^N$, that is,
\be\label{eq:weakharm}
\int_\Omega \sum_{k=1}^M (u_i \nabla_k u_j-u_j\nabla_k u_i)\nabla_k\varphi^{ij}\,dx=0, \quad  \forall \varphi^{ij}\in C^\infty_c(\Omega;\R), 1\leq i<j\leq N.
\ee
To prove \eqref{eq:weakharm}, by multiplying \eqref{eq:eps} by $u_{\eps,j} \varphi^{ij}\in H^1_0(\Omega)$ and integrating by parts, we have
$$
\int_\Omega \nabla u_{\eps,i} \cdot \nabla (u_{\eps,j}\varphi^{ij})\,dx=\frac1{\eps^2} \int_\Omega W'(1-|u_\eps|^2) u_{\eps,i}  \cdot u_{\eps,j}\varphi^{ij}\, dx,  \quad \forall 1\leq i<j\leq N.
$$ Setting $\tilde\varphi^{ji}:=\varphi^{ij}$, by reversing above $i$ and $j$ and then by subtraction we obtain 
$$
\int_\Omega \sum_{k=1}^M (u_{\eps,i} \nabla_k u_{\eps,j}-u_{\eps,j}\nabla_k u_{\eps,i})\, \nabla_k\varphi^{ij}\,dx=0,  \quad \forall 1\leq i<j\leq N.
$$ As $u_\eps\rightharpoonup u$ weakly in $H^1$ and strongly in $L^2$, passing to the limit $\eps\to 0$, we get  \eqref{eq:weakharm}. 

Assume now that \eqref{ass:convH1} and \eqref{ass:Wsweak} hold. Then
\be\label{eq:3456}
\int_\Omega |\nabla u|^2\, dx\leq \liminf_{\eps\to 0} \int_\Omega |\nabla u_\eps|^2\, dx \leq \limsup_{\eps\to 0} \int_\Omega |\nabla u_\eps|^2\, dx\leq 
2 \limsup_{\eps\to 0} \mcF_\eps[u_\eps] = \int_\Omega |\nabla u|^2\, dx.
\ee
Hence $\nabla u_\eps\to \nabla u$ in $L^2(\Omega)$. Returning to \eqref{ass:convH1}, we deduce that $\frac{1}{\vareps^2}W(1 - |u_\vareps|^2) \rightarrow 0$ in $L^1(\Omega)$. 
It remains to show that $S$ is divergence-free. Consider the stress-energy tensor $S^\vareps$ of $u_\vareps$ associated to the functional $\mcF_\vareps$ defined for $i,j=1, \dots, M$:
\[
S^\vareps_{ij} = e_\vareps(u_\vareps)\delta_{ij} - \nabla_i u_\vareps \cdot \nabla_j u_\vareps,
\]
where $e_\vareps$ is defined in \eqref{def:edensity}.
By \eqref{ass:convH1}, $e_\vareps(u_\vareps) \rightarrow \frac{1}{2}|\nabla u|^2$ in $L^1(\Omega)$ and so $S^\vareps \rightarrow S$ in $L^1(\Omega)$ as $\eps\to 0$. Therefore, it is enough to show that $S^\eps$ is divergence-free. Since $W \in C^1$ and $u_\vareps \in H^1 \cap L^\infty(\Omega)$, by elliptic regularity for \eqref{eq:eps}, we have $u_\eps \in W^{2,p}_{\rm loc}(\Omega)$ for any $p \in [1,\infty)$ and \eqref{eq:eps} holds almost everywhere. This also implies that $S^\vareps \in W^{1,p}_{\rm loc}(\Omega)$ for any $p \in [1,\infty)$. Using \eqref{eq:eps}, we compute for $1\leq j \leq N$:
\begin{align*}
\sum_{i=1}^M \nabla_i S^\vareps_{ij}
	&=  \sum_{i=1}^M\nabla_i \nabla_j u_\vareps \cdot \nabla_i u_\vareps - \frac{1}{\vareps^2} W'(1 - |u_\vareps|^2) u_\vareps \cdot  \nabla_j u_\vareps \\
		&\qquad -  \sum_{i=1}^M\big(\nabla_i \nabla_i u_\vareps \cdot \nabla_j u_\vareps + \nabla_i u_\vareps \cdot \nabla_i \nabla_j u_\vareps\big)\\
	&=  -\Big(\Delta u_\vareps + \frac{1}{\vareps^2} W'(1 - |u_\vareps|^2) u_\vareps\Big) \cdot  \nabla_j u_\vareps
	\stackrel{\eqref{eq:eps}}{=} 0 \quad \textrm{a.e. in $\Omega$}.
\end{align*}
The proof is complete.
\end{proof}

\begin{lemma}
\label{Lem:u0ConvSHM}
Let $\Omega \subset \mathbb{R}^M$ be a $C^2$ bounded domain,  $W \in C^2((-\infty,1])$ satisfy \eqref{ass:Wstructure}, $(u_\eps)_{\eps\downarrow 0}$ be critical points of $\mcF_\eps$ in $H^1\cap L^\infty(\Omega,\RR^N)$ satisfying \eqref{ass:convH1} for some $u_0 \in H^1(\Omega,\RR^N)$. Suppose also that $u_\eps|_{\partial\Omega} \in C^0(\partial\Omega, \mathbb{S}^{N-1})$. Then $u_\eps \rightarrow u_0$ in $H^1(\Omega,\RR^N)$, $|u_0| = 1$ a.e. in $\Omega$, and $u_0$ is a stationary $\mathbb{S}^{N-1}$-valued harmonic map.
\end{lemma}

\begin{proof}
By inspection of \eqref{eq:3456} (which does not require that $|u_0|=1$ in $\Omega$), we have $u_\eps \rightarrow u_0$ in $H^1(\Omega)$ and 
$\frac{1}{\vareps^2}W(1 - |u_\vareps|^2) \rightarrow 0$ in $L^1(\Omega)$. In particular, for a sequence $\eps\to 0$, $u_\eps\to u_0$ a.e. in $\Omega$ and $\frac{1}{\vareps^2}W(1 - |u_\vareps|^2) \to 0$ a.e. in $\Omega$. Thus, $W(1-|u_0|^2)=0$ a.e. in $\Omega$.
Recall that $|u_\eps| \leq 1$ by Theorem \ref{thm:MPI}; thus, $|u_0|\leq 1$ a.e. in $\Omega$. As $W > 0$ in $(0,1]$ (by \eqref{ass:Wstructure}), this implies  $|u_0|=1$ a.e. in $\Omega$. Finally, by Lemma \ref{lem:Stationarity},  we conclude that $u_0$ is a stationary $\mathbb{S}^{N-1}$-valued harmonic map.
\end{proof}

\subsection{An elliptic boundary estimate}

In the following we provide an elliptic boundary estimate in the half ball 
$$B_r^+(0)=\{x\in \R^M\, :\, |x|<r, x_M>0\}$$ which is a slightly changed version of \cite[Lemma $11$]{NZ13}. We denote $x=(x',x_M)\in \R^M$ with $x'\in \R^{M-1}$.

\begin{lemma}
\label{lem:bdry-est}
Let $\beta \in (0,1)$, $(a_{ij})_{1 \leq i,j \leq M} \in C^{0,\beta}(B_r^+(0))$ be uniformly elliptic, $(b_i)_{1 \leq i \leq M} \in L^\infty(B_r^+(0))$ and let $u \in H^{1}(B_r^+(0))$ be a weak solution to
\[
\begin{cases}
- \sum_{i,j = 1}^M \nabla_i(a_{ij} \nabla_j u) + \sum_{i=1}^M b_i \nabla_i u = f + \sum_{i=1}^M \nabla_i h_i \text{ in } B_r^+(0),\\
u = g \text{ on } B_r(0) \cap\{x_M = 0\},
\end{cases}
\]
where $f \in L^p(B_r^+(0))$ with $p > M$, $h_i \in C^{0,\beta}(B_r^+(0))$ and $g \in C^{1,\beta}( B_r(0) \cap\{x_M = 0\})$. Then there exists a constant $C$ depending only on $M$, $\beta$, $p$, the ellipticity constants of $(a_{ij})$, and an upper bound for $\|a_{ij}\|_{C^{0,\beta}(B_r^+(0))}$ and $\|b_{i}\|_{L^\infty(B_r^+(0))}$ such that
\begin{align*}
\|\nabla u\|_{C^{0}(\bar B_{r/2}^+(0))} + r^\beta [\nabla u]_{C^{0,\beta}(B_{r/2}^+(0))} 
	&\leq C\Big[r^{1 - \frac{M}{p}} \|f\|_{L^p(B_r^+(0))} 
		+ \sum_{i=1}^M \Big(\|h_i\|_{C^{0}(\bar B_{r }^+(0))} + r^\beta [h_i]_{C^{0,\beta}(B_{r }^+(0))}\Big)\\
		&\quad 
		+ \|\nabla_\tau g\|_{C^0(\bar B_r(0) \cap\{x_M = 0\})} 
		+ r^\beta[\nabla_\tau g]_{C^{0,\beta}(B_r(0) \cap\{x_M = 0\})}\\
		&\qquad
		+ r^{-\frac{M}{2}} \|\nabla u\|_{L^2(B_r^+(0))} \Big].
\end{align*}
\end{lemma}

\begin{proof}
By scaling, it suffices to consider $r = 1$. Also, replacing $u$ by $u(x', x_M) - g(x',0)$, $f$ by $f(x',x_M) - \sum_{i=1}^{M-1} b_i(x',x_M) \nabla_i g(x',0)$ and $h_i$ by $h_i(x',x_M) + \sum_{j=1}^{M-1} a_{ij}(x',x_M)\nabla_j g(x',0)$, we may assume without loss of generality that $g \equiv 0$. Then we have the chain of inequalities:
\begin{align*}
\|\nabla u\|_{C^{0,\beta}(B_{1/2}^+(0))} 
	&\leq C\Big[  \|f\|_{L^p(B_{3/4}^+(0))} 
		+ \sum_{i=1}^M  \|h_i\|_{C^{0,\beta}(B_{3/4}^+(0))}  +   \| u\|_{L^\infty(B_{3/4}^+(0))} \Big],\\
\|u\|_{L^\infty(B_{3/4}^+(0))} 
	&\leq C\Big[  \|f\|_{L^p(B_{1}^+(0))} 
		+ \sum_{i=1}^M  \|h_i\|_{C^{0}(\bar B_{1}^+(0))}  +   \| u\|_{L^2(B_1^+(0))} \Big],\\
\| u\|_{L^2(B_1^+(0))}
	&\leq C\|\nabla u\|_{L^2(B_1^+(0))}.	
\end{align*}
where we have used \cite[Corollary 8.36]{GT} for the first inequality, \cite[Theorem 8.25]{GT} for the second inequality, and Poincar\'e's inequality for the third inequality. The conclusion follows.
\end{proof}

 \def\cprime{$'$}

\end{document}